\documentclass{siamltex}

\usepackage{amsmath}
\usepackage{amsfonts}
\usepackage{amssymb}

\usepackage{tikz}
\usepackage{pgfplots}
\usetikzlibrary{plotmarks,calc}

\newcommand{\mb}\boldsymbol

\newcommand{\pdiff}[2]{\frac{\partial #1}{\partial #2}}
\newcommand{\pddiff}[2]{\frac{\partial^2 #1}{\partial #2^2}}

\newcommand{\ie}{\emph{i.e.}~}
\newcommand{\eg}{\emph{e.g.}~}
\newcommand{\etc}{\emph{etc.}}

\newcommand{\Ord}[1]{\mathcal{O}\left(#1\right)}

\newcommand{\tOrd}[1]{\mathcal{O}\big(#1\big)}

\newcommand{\signum}{\operatorname{\text{sgn}}}

\def\Xint#1{\mathchoice
   {\XXint\displaystyle\textstyle{#1}}%
   {\XXint\textstyle\scriptstyle{#1}}%
   {\XXint\scriptstyle\scriptscriptstyle{#1}}%
   {\XXint\scriptscriptstyle\scriptscriptstyle{#1}}%
   \!\int}
\def\XXint#1#2#3{{\setbox0=\hbox{$#1{#2#3}{\int}$}
     \vcenter{\hbox{$#2#3$}}\kern-.5\wd0}}

\def\dashint{\Xint-}

\usepackage{graphicx}

\usepackage{color}
\usepackage[normalem]{ulem}

\definecolor{darkgreen}{RGB}{0,127,0}

\newcommand{\rout}[1]{}

\usepackage{hyperref}
\hypersetup{
    unicode=false,          
    pdftoolbar=true,        
    pdfmenubar=true,        
    pdffitwindow=true,      
    pdftitle={},    
    pdfauthor={Cameron Hall},     
    pdfsubject={},   
    pdfnewwindow=true,      
    pdfkeywords={}, 
    colorlinks=true,       
    linkcolor=black,          
    citecolor=black,        
    filecolor=black,      
    urlcolor=black,           
}

\title{The mechanics of a chain or ring of spherical magnets}

\author{Cameron L.~Hall, Dominic Vella and Alain Goriely \\ Accepted for Publication in SIAM Journal on Applied Mathematics}

\date{\today}

\begin{document}

\maketitle

\begin{abstract}
Strong magnets, such as neodymium-iron-boron magnets, are increasingly being  manufactured as spheres. Because of their dipolar characters, these spheres can easily be arranged into long chains that exhibit mechanical properties reminiscent of  elastic strings or rods. While simple formulations exist for the energy of a deformed elastic rod, it is not clear whether or not they are also appropriate for a chain of spherical magnets. In this paper, we use discrete-to-continuum asymptotic analysis to derive a continuum model for the energy of a deformed chain of magnets based on the magnetostatic interactions between individual spheres. We find that the mechanical properties of a chain of magnets differ significantly from those of an elastic rod: while both magnetic chains and elastic rods support bending by change of local curvature, nonlocal interaction terms also appear in the energy formulation for a magnetic  chain. This continuum model for the energy of a chain of magnets is used to analyse small 
deformations of a circular ring of magnets and hence obtain theoretical predictions for the vibrational modes of a circular ring of magnets. Surprisingly, despite the contribution of nonlocal energy terms, we find that the vibrations of a circular ring of magnets are governed by the same equation that governs the vibrations of a circular elastic ring.
\end{abstract}



\section{Introduction}
\label{S:Intro}


From hard drives to electric motors, neodymium-iron-boron (NdFeB) magnets have become ubiquitous in low-temperature applications where high magnetic strength is required: in 2008, the global production of NdFeB magnets exceeded 60 000 tons \cite{Gutfleisch2011}. In the last few years, collections of spherical NdFeB magnets have been sold as toys under brand names including Neocube\texttrademark and Buckyballs\texttrademark. The high magnetic strength and low weight of these spherical magnets means that they can be used to construct complicated and interesting structures that are held together by the magnetic attraction between spheres. 

\begin{figure}
 \centering
 \includegraphics[width=0.6\textwidth]{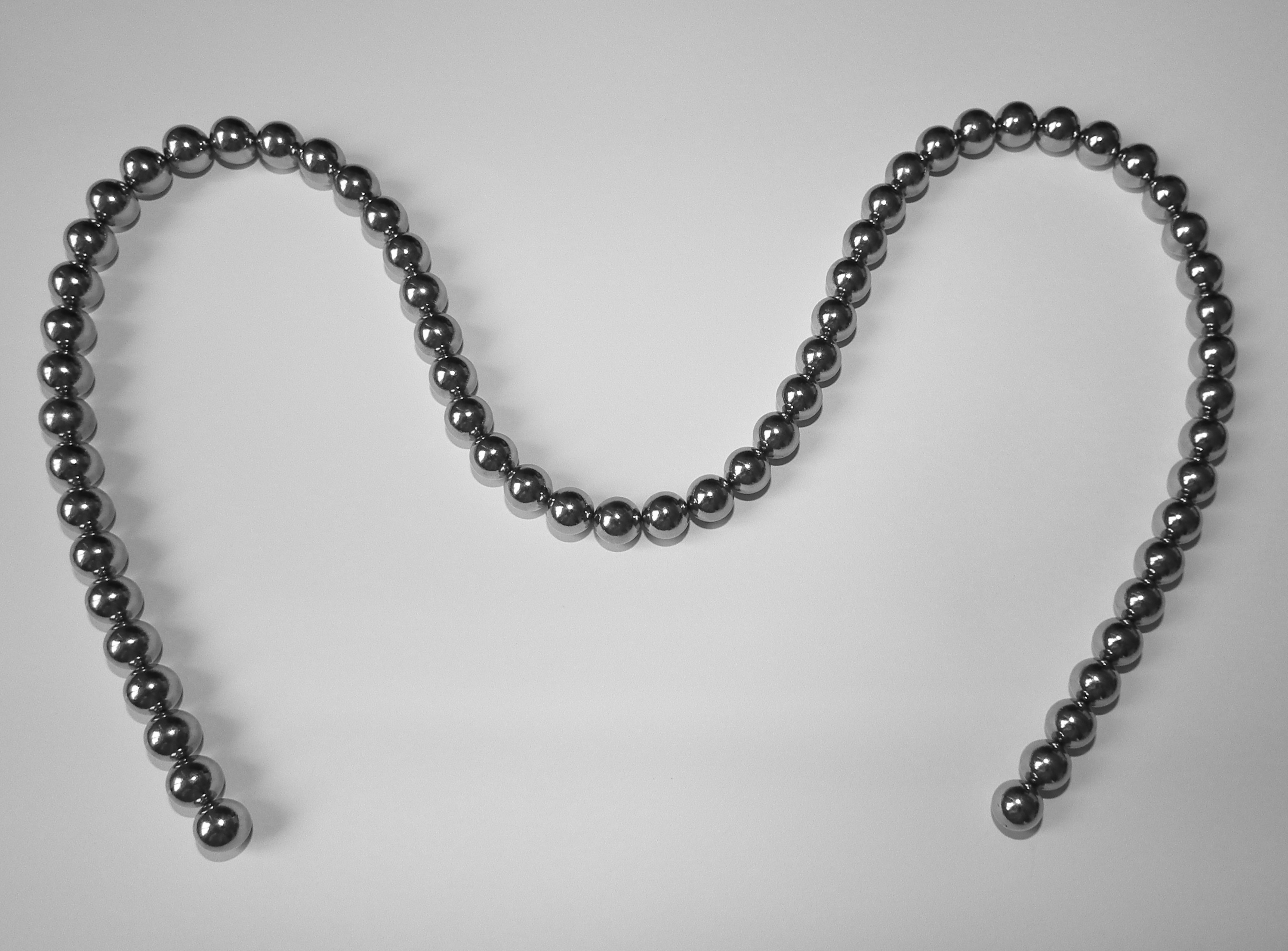}
 \caption{Photograph of a chain of spherical magnets resting on a table-top. The diameter of the spheres in this image is $5\mathrm{~mm}$. (Image courtesy of Emmanuel du Pontavice.)}
 \label{F:ChainPhoto}
\end{figure}

Perhaps the simplest structure that can be built from NdFeB spherical magnets is a chain, as shown in Figure \ref{F:ChainPhoto}. While such chains are clearly composed of discrete particles (the individual magnets), they exhibit mechanical behaviour that is reminiscent of elastic rods. For example, straight chains of magnets buckle reversibly if a sufficient load is applied, while a `squashed' circular ring of magnets will undergo damped vibrations and return to a circular shape if the applied load is released. This naturally leads to a number of questions: What is an appropriate mechanical model for a chain of spherical magnets? Do chains of magnets have an effective, magnetically induced, `bending stiffness'? How do the mechanical properties of a chain depend on the field strength of the magnets? In this paper, we address these questions by using a discrete-to-continuum asymptotic analysis 
to develop a continuum model of a chain of spherical magnets based on the magnetostatic interactions between individual spheres.

Spherical magnets with a uniform magnetization produce an external field that is precisely that of a dipole  \cite{Jackson}. The interactions between dipolar hard spheres have been well-studied using both numerical and physical simulations. For example, \cite{Clarke1994} uses numerical methods to determine the ground-state equilibria for systems containing small numbers of dipolar hard spheres, while \cite{Guo2005} studies the ring sizes obtained after experimentally vibrating a large collection of magnetic balls and \cite{Ku2011} uses both experiments and simulations to explore the self-assembly of magnetic particles into crystalline structures. In all of these papers, the focus is on discrete-scale models, where the position and orientation of each dipolar particle are tracked individually.

An alternative approach is to construct continuum models based on the interactions of many magnets. For example, \cite{Yoon2010} presents continuum models for aggregates of magnetic dipoles; however, the authors do not relate the continuum properties of their system (\eg the bending stiffness of a chain of spherical magnets), to the properties of the individual magnets (\eg their radius and magnetic field strength). An advance on this can be achieved by obtaining a `bending stiffness' from a comparison of the energy of a single magnet in an infinitely-long straight chain with the energy of a single magnet in a circular ring of specified radius \cite{VellaPreparation}.  Unfortunately, this calculation says nothing about whether a bending stiffness formulation is appropriate for a chain with varying radius of curvature. 

In this paper, we use asymptotic methods to derive a continuum model of a chain of spherical magnets based on discrete-scale dipole interactions. This means that we do not need to make  \emph{a priori} assumptions about the form of our model (\eg that energy is everywhere proportional to the square of the curvature of the chain), and we ultimately obtain expressions for the energy of a chain that are not intuitively obvious. Thus, our asymptotic approach enables us to gain insights into the physics of a chain of magnets that could not be easily obtained using other methods. 

Our results depend on an analysis of the sum that represents the energy of a chain of spherical magnets. By assuming that the radius of curvature of the chain is large relative to the radius of an individual magnet, and by assuming that the orientation of the magnetic dipoles varies slowly along the chain, we find that we can approximate this sum with a continuous functional. 
Then, the problem of finding the locations of the individual magnets is reduced to a variational problem for a continuous curve representing the shape of the chain.

There is a strong parallel between this approach and the work described by Hall and coworkers \cite{Hall2011,Hall2010}, in which asymptotic methods are used to describe the interactions of a large number of repelling particles (such as dislocations or dislocation dipoles) in a `pile-up' configuration. As in the present work, \cite{Hall2011} and \cite{Hall2010} involves approximating a sum 
by exploiting assumptions about the positions of the particles. 
This ultimately led to a continuum model for the density of particles within the system.

In \cite{Hall2010}, it was found that the form of the continuum model depends strongly on the rate at which the inter-particle repulsion decays with distance.  Since spherical magnets can be modelled as magnetic dipoles and the force of attraction between dipoles decays in inverse proportion to the fourth power of distance (see, for example, \cite{Jackson}), the results in \cite{Hall2010} might lead us to expect local interactions within a chain of magnets to dominate over nonlocal interactions. Thus, the main contribution to the energy of a chain of spherical magnets should be analogous to the bending energy of a rod or beam, or at least depend only on the curvature and higher derivatives of the shape of the chain. 

However, an unexpected result of our work is that both local and nonlocal contributions to the energy are important. As described in Section \ref{S:Main}, a chain of spherical magnets does have an effective bending stiffness arising from local dipole interactions, but the energy associated with this resistance to 
local bending appears at the same asymptotic order as a nonlocal energy associated with long-range dipole interactions. 

Another observation in \cite{Hall2011,Hall2010} is that continuum dislocation models are inappropriate for describing the ends of a dislocation pile-up, and it is instead necessary to consider discrete problems in the boundary layers. Similarly, we find that the continuum model developed in this paper for a chain of magnets breaks down near the ends of the chain. In the present work, we restrict our analysis to cases where a continuum model is valid. While this prevents us from developing a complete model of a finite chain of spherical magnets, it is possible for us to perform a full analysis of a ring of spherical magnets. In Section \ref{S:Rings},  we calculate the energy of a circular ring of magnets and we consider the effects of a small perturbation to the circular shape. This enables us to obtain equations for the in-plane vibrational modes of a circular ring of magnets, and we recover the nonintuitive result 
that these are identical to the in-plane vibrational modes of an elastic ring. Our results match with the theoretical and experimental results obtained in \cite{VellaPreparation}.

\section{Problem statement}
\label{S:ProblemStatement}

\subsection{Physical background}
\label{S:PhysicalBackground}

Consider a chain of $n + 1$ identical spherical magnets with centres located at $\{\mb{r}_i\}$ where $i = 0,\,1,\,2,\,\ldots,\,n$. As shown in \cite{Jackson} (p.~198) and other standard texts on electromagnetism, a uniformly magnetized sphere of radius $a$ and magnetization intensity $\mb{I}$ generates an exterior magnetic field identical to that generated by a magnetic dipole with moment $\mb{m} = \frac{4}{3} \, \pi \, a^3 \, \mb{I}$. Typically, spherical NdFeB magnets are sold according to their characteristic magnetic field strength, $B = \mu_0 \, |\mb{I}|$, where $\mu_0$ is the magnetic permeability of free space. Hence, we find that the magnetic field of the $i$th magnet is identical to that of a dipole with moment
\begin{equation}
  \mb{m}_i = \frac{4 \, \pi \, a^3 \, B}{3 \, \mu_0} \, \mb{\hat{m}}_i,
\end{equation}
where $\mb{\hat{m}}_i$ is a unit vector in the direction of the dipole of the $i$th magnet.

The exterior field due to a dipole in free space located at $\mb{r}_i$ with moment $\mb{m}_i$ is given by
\begin{equation}
 \mb{B}_i(\mb{r}) = \frac{\mu_0}{4 \, \pi} \, \frac{3 \, \big[(\mb{r} - \mb{r}_i) \cdot \mb{m}_i\big] \, (\mb{r} - \mb{r}_i) - ||\mb{r} - \mb{r}_i||^2 \, \mb{m}_i}{||\mb{r} - \mb{r}_i||^5}. \label{ExtField_Dim}
\end{equation}
Since the relative permeability of NdFeB magnets is close to unity (around 1.05 in most manufacturers' specifications), this exterior field will not be significantly altered by the presence of more magnets. Hence, the fields from several magnets can be summed to give the overall magnetic field. In the absence of an  applied magnetic field, it follows that the total magnetic field at a point $\mb{r}$ due to a chain of $n$ magnets with centres at $\{\mb{r}_i\}$ and dipole moments of $\{\mb{m}_i\}$ is given by
\begin{equation}
 \mb{B}_\text{tot}(\mb{r}) 
 = \sum_{i=0}^n \mb{B}_i(\mb{r}) 
 = \frac{\mu_0}{4 \, \pi} \, \sum_{i=0}^n \frac{3 \, \big[(\mb{r} - \mb{r}_i) \cdot \mb{m}_i\big] \, (\mb{r} - \mb{r}_i) - ||\mb{r} - \mb{r}_i||^2 \, \mb{m}_i}{||\mb{r} - \mb{r}_i||^5}.
\end{equation} 

As described in \cite{Jackson} (p.~190), the energy of a dipole with moment $\mb{m}$ introduced into a magnetic field $\mb{B}_0$ is given by
\begin{equation}
 E = - \mb{m} \cdot \mb{B}_0.\label{DipoleEnergy}
\end{equation} 
It follows that the energy of a body with magnetization intensity $\mb{I}$ introduced into a field $\mb{B}_0$ is
\begin{equation}
 E = - \iiint \mb{I} \cdot \mb{B}_0 \, dV, 
\end{equation}
where the integral is taken over the volume of the magnetized body. If a uniformly magnetized sphere of radius $a$ and magnetization intensity $\mb{I} = \left(\frac{4}{3} \, \pi \, a^3\right)^{-1} \, \mb{m}$ is introduced into the field due to a dipole of strength $\mb{m}_0$ located at the origin, its energy will be given by
\begin{equation}
 E = - \mb{m}^T \, \frac{3}{4 \, \pi \, a^3} \iiint \frac{3 \, \mb{r} \otimes \mb{r} - ||\mb{r}||^2 \, \mb{\delta}}{||\mb{r}||^5} \, dV \, \mb{m}_0,
\end{equation}
where $\mb{\delta}$ is the Kronecker delta tensor. This integral can be evaluated analytically, yielding the result that
\begin{equation}
 E = -\mb{m}^T \, \frac{3 \, \mb{R} \otimes \mb{R} - ||\mb{R}||^2 \, \mb{\delta}}{||\mb{R}||^5} \, \mb{m}_0,
\end{equation}
where $\mb{R}$ is the position vector of the centre of the sphere. Thus, by comparison with  \eqref{ExtField_Dim} and \eqref{DipoleEnergy}, we find that the energy of interaction between two spherical magnets is identical to the energy of interaction between two dipoles.

Noting that we need to avoid double-counting the energy of interaction between two magnets, we can therefore use \eqref{DipoleEnergy} to obtain the total energy of a chain as follows:
\begin{equation}
 E_\text{tot} 
 = \frac{1}{2} \, \sum_{i=0}^n E_i 
 = -\frac{1}{2} \, \sum_{i=0}^n \mb{m}_i \cdot \text{reg} \left[ \mb{B}_\text{tot}(\mb{r}) \right]_{\mb{r} = \mb{r}_i}
 = -\sum_{i=0}^n \, \sum_{\substack{j=0 \\ j\neq i}}^{n} \frac{\mb{m}_i \cdot \mb{B}_j(\mb{r}_i)}{2},
 \label{Etot_Temp}
\end{equation} 
where $\text{reg} \left[ \mb{B}_\text{tot}(\mb{r}) \right]$ represents the regular part of $\mb{B}_\text{tot}$, defined so that
\begin{equation}
 \text{reg} \left[ \mb{B}_\text{tot}(\mb{r}) \right]_{\mb{r} = \mb{r}_i} = \lim_{\mb{r} \rightarrow \mb{r}_i} \left[\mb{B}_\text{tot}(\mb{r}) - \mb{B}_i(\mb{r}) \right].
\end{equation}
This is necessary because $\mb{B}_\text{tot}(\mb{r})$ is singular at each $\mb{r} = \mb{r}_i$, and because the energy of each magnet depends on its interactions with the external magnetic field, not with its own magnetic field.

Throughout this paper, we ignore the effects of electromagnetic induction and assume that the potential energy of the entire system is given simply by equation \eqref{Etot_Temp}. While this is correct for static systems, we also consider the dynamics of an oscillating ring in Section \ref{S:Rings}. However, we can justify continuing to neglect induction by noting that the dissipation of energy through the small current induced in the magnets will be negligible compared to the magnetic potential energy given by \eqref{Etot_Temp} and the kinetic energy based on the mass and velocity of the magnets. 

A further issue is the problem of contact friction between spherical magnets. Simple experiments with marked magnets strongly indicate that magnets initially in contact will slide and roll against each other in order to minimise the energy of the system, regardless of frictional forces that might resist movement. While friction has the potential to cause significant damping in dynamical problems, including the oscillating ring described in Section \ref{S:Rings}, we ignore its effects in the present analysis.

\subsection{Continuum formulation}
\label{S:ContinuumFormulation}

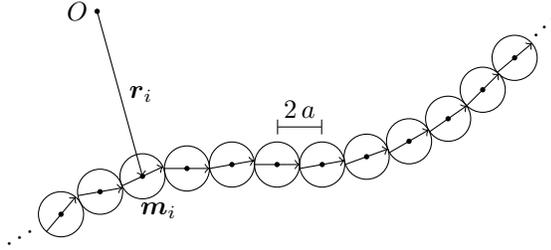
\begin{figure}
 \centering
 \begin{tikzpicture}
  \coordinate[label=left:{$O$}] (origin) at (0,0);
  \node[rotate=30] (R0m) at (-1,-3) {$\cdots$};
  \draw ($(R0m)+(30:0.6cm)$) coordinate (R1) circle (0.3cm);
  \draw [->] ($(R1)-(50:0.3cm)$) -- ($(R1)+(50:0.3cm)$); 
  \draw ($(R1)+(30:0.6cm)$) coordinate (R2) circle (0.3cm);
  \draw [->] ($(R2)-(10:0.3cm)$) -- ($(R2)+(10:0.3cm)$); 
  \draw ($(R2)+(20:0.6cm)$) coordinate (R3) circle (0.3cm);
  \draw [->] ($(R3)-(25:0.3cm)$) -- ($(R3)+(25:0.3cm)$); 
  \node at ($(R3)-(115:0.5cm)$) {$\mb{m}_i$};
  \draw ($(R3)+(10:0.6cm)$) coordinate (R4) circle (0.3cm);
  \draw [->] ($(R4)-(0:0.3cm)$) -- ($(R4)+(0:0.3cm)$); 
  \draw ($(R4)+(5:0.6cm)$) coordinate (R5) circle (0.3cm);
  \draw [->] ($(R5)-(10:0.3cm)$) -- ($(R5)+(10:0.3cm)$); 
  \draw ($(R5)+(0:0.6cm)$) coordinate (R6) circle (0.3cm);
  \draw [->] ($(R6)-(0:0.3cm)$) -- ($(R6)+(0:0.3cm)$); 
  \draw ($(R6)+(0:0.6cm)$) coordinate (R7) circle (0.3cm);
  \draw [->] ($(R7)-(10:0.3cm)$) -- ($(R7)+(10:0.3cm)$); 
  \draw ($(R7)+(10:0.6cm)$) coordinate (R8) circle (0.3cm);
  \draw [->] ($(R8)-(20:0.3cm)$) -- ($(R8)+(20:0.3cm)$); 
  \draw ($(R8)+(20:0.6cm)$) coordinate (R9) circle (0.3cm);
  \draw [->] ($(R9)-(30:0.3cm)$) -- ($(R9)+(30:0.3cm)$); 
  \draw ($(R9)+(30:0.6cm)$) coordinate (R10) circle (0.3cm);
  \draw [->] ($(R10)-(35:0.3cm)$) -- ($(R10)+(35:0.3cm)$); 
  \draw ($(R10)+(40:0.6cm)$) coordinate (R11) circle (0.3cm);
  \draw [->] ($(R11)-(45:0.3cm)$) -- ($(R11)+(45:0.3cm)$); 
  \draw ($(R11)+(45:0.6cm)$) coordinate (R12) circle (0.3cm);
  \draw [->] ($(R12)-(40:0.3cm)$) -- ($(R12)+(40:0.3cm)$); 
  \node[rotate=45] (Rend) at ($(R12)+(45:0.6cm)$) {$\cdots$};
  \foreach \point in {origin,R1,R2,R3,R4,R5,R6,R7,R8,R9,R10,R11,R12}
    \fill [black] (\point) circle (1pt);
  \draw [->] (origin) to node [right]{$\mb{r}_i$} (R3);
  \draw [|-|] ($(R6)+(90:0.5cm)$) to node[above]{$2 \, a$} ($(R7)+(90:0.5cm)$);
 \end{tikzpicture}
 \caption{Schematic diagram of a chain of spherical magnets. The $i$th magnet in the chain has a position vector of $\mb{r}_i$ and a dipole vector of $\mb{m}_i$. The constraint that the spheres form a chain ensures that the distance between the centres of any two neighbouring magnets is given by $2 \, a$.}
 \label{F:ChainSchematic}
\end{figure}

Since adjacent spherical magnets in the chain will be in direct contact with each other, the centre points will be regularly separated so that $||\mb{r}_{i\pm1} - \mb{r}_i|| = 2 \, a$. This configuration is illustrated in Figure \ref{F:ChainSchematic}. 

We are interested in investigating the case where the chain contains a large number of magnets and where the shape of the chain can be approximated by a smooth curve. More formally, we make the fundamental assumption that $\mb{r}_i$ and $\mb{m}_i$ can be expressed in the form
\begin{equation}
 \mb{r}_i = 2 \, a \, n \, \mb{\tilde{r}}(i\, n^{-1}; \, n), \quad  \mb{m}_i = \frac{4 \, \pi \, a^3 \, B}{3 \, \mu_0} \, \mb{\tilde{m}}(i\,n^{-1}; \, n), \label{NondimMain}
\end{equation}
where $\mb{\tilde{r}}(s;\,n)$ and $\mb{\tilde{m}}(s;\,n)$ are suitably differentiable dimensionless functions that are expressed as asymptotic expansions as $n \rightarrow \infty$, and $2 \, a \, n$ represents the length of the entire chain (defined by the sum of straight segments connecting the centres of the magnets). It follows from these definitions that 
\begin{equation}
||\mb{\tilde{r}}((i\pm1) \, n^{-1};\,n) - \mb{\tilde{r}}(i \, n^{-1};\,n)|| = n^{-1}, \label{NondimGapNeighbours}
\end{equation}
and
\begin{equation}
 ||\mb{\tilde{m}}(s;\,n)|| = 1. \label{mNorm=1}
\end{equation}


In general, we nondimensionalize distances with the chain length $2\,a\,n$ so that $\mb{r} = 2 \, a \, n \, \mb{\bar{r}}$, and we nondimensionalize magnetic moments with their magnitude $ \frac{4\pi}{3} \, a^3 \, B/\mu_0$. Based on these definitions, we find that appropriate nondimensionalizations for $\mb{B}_\text{tot}$ and $E_\text{tot}$ are given by
\begin{equation}
 \mb{B}_\text{tot}(\mb{r}) = \frac{B}{24} \, \mb{\bar{B}}_\text{tot}(\bar{\mb{r}}), \quad E_\text{tot} = \frac{\pi \, a^3 \, B^2}{18 \, \mu_0} \bar{E}_\text{tot}, 
\end{equation}
where $\mb{\bar{B}}_\text{tot}(\mb{\bar{r}})$ and $\bar{E}_\text{tot}$ are dimensionless representations of the total magnetic field and total energy respectively.

With this notation, we obtain
\begin{multline}
 \mb{\bar{B}}_\text{tot}(\mb{\bar{r}}) 
 = n^{-3} \, \sum_{i=0}^n \frac{3 \, \big[(\mb{\bar{r}} - \mb{\tilde{r}}(i\,n^{-1};\,n)) \cdot \mb{\tilde{m}}(i\,n^{-1}; \, n)\big] \, (\mb{\bar{r}} - \mb{\tilde{r}}(i\,n^{-1};\,n))}{||\mb{\bar{r}} - \mb{\tilde{r}}(i\,n^{-1};\,n)||^5} \\ 
 - n^{-3} \, \sum_{i=0}^n \frac{\mb{\tilde{m}}(i\,n^{-1}; \, n)}{||\mb{\bar{r}} - \mb{\tilde{r}}(i\,n^{-1};\,n)||^3}, \label{BtotNondimOrig}
\end{multline}
and
\begin{equation}
 \bar{E}_\text{tot} = - \frac{1}{2} \, \sum_{i=0}^n \mb{\tilde{m}}_i \cdot \text{reg} \left[ \mb{\bar{B}}_\text{tot}(\mb{\bar{r}}) \right]_{\mb{\bar{r}} = \mb{\tilde{r}}(i\,n^{-1};\,n)} \label{EtotNondim}
\end{equation}

Lastly, it is useful  to introduce the functions $\mb{\tilde{B}}(s;\,n)$ and $\tilde{E}(s;\,n)$ to represent the magnetic field and energy density, respectively, along the centre line of the chain of magnets. These are defined so that
\begin{multline}
 \mb{\tilde{B}}(s;\,n) = \sum_{i=0}^n \mb{\tilde{B}}_i(s)\\
 = n^{-3} \, \sum_{i=0}^n \frac{3 \, \big[(\mb{\tilde{r}}(s;\, n) - \mb{\tilde{r}}(i\,n^{-1};\,n)) \cdot \mb{\tilde{m}}(i\,n^{-1}; \, n)\big] \, (\mb{\tilde{r}}(s;\, n) - \mb{\tilde{r}}(i\,n^{-1};\,n))}{||\mb{\tilde{r}}(s;\, n) - \mb{\tilde{r}}(i\,n^{-1};\,n)||^5} \\ 
 - n^{-3} \, \sum_{i=0}^n \frac{\mb{\tilde{m}}(i\,n^{-1}; \, n)}{||\mb{\tilde{r}}(s;\, n) - \mb{\tilde{r}}(i\,n^{-1};\,n)||^3}, \label{BNondimLong}
\end{multline}
and
\begin{equation}
 \tilde{E}(s;\,n) = - \mb{m}(s;\,n) \cdot \mb{\tilde{B}}_\text{reg}(s;\,n), \label{EDefn}
\end{equation} 
where $\mb{\tilde{B}}_\text{reg}(s;\,n)$ is an appropriately smooth and slowly varying function with the property that
\begin{equation}
 \mb{\tilde{B}}_\text{reg}(i \, n^{-1};\,n) 
 = \lim_{s \rightarrow i \, n^{-1}} \left[ \mb{\tilde{B}}(s;\,n)  - \mb{\tilde{B}}_i(s;\,n) \right], \label{BRegDefn}
\end{equation}
and $\mb{\tilde{B}}_i(s;\,n)$ is the summand in \eqref{BNondimLong}, obtained by substituting $\mb{r} = \mb{\tilde{r}}(s;\,n)$ into \eqref{ExtField_Dim} and nondimensionalizing.

Throughout this section, we have used the notation $\mb{\tilde{r}} = \mb{\tilde{r}}(s;\,n)$, $\mb{\tilde{m}} = \mb{\tilde{m}}(s;\,n)$ \etc~to emphasize the fact that $\mb{\tilde{r}}$, $\mb{\tilde{m}}$, $\mb{\tilde{B}}$ and ${\tilde{E}}$ are all dependent on $n$. In the next section, these functions will be expressed as asymptotic expansions for large $n$; for conciseness of notation, we will omit the explicit dependence on $n$ from this point onwards.

While our ultimate aim is to obtain an expression for $\bar{E}_\text{tot}$, the first step in our analysis will be the construction of an asymptotic representation of $\mb{\tilde{B}}(s)$. As the definition of $\mb{\tilde{B}}(s)$ in \eqref{BNondimLong} involves a simple sum over all values of $i$, it is easier to analyse $\mb{\tilde{B}}(s)$ than to analyse directly the double-sum with excluded terms that defines the total energy in \eqref{Etot_Temp}. Moreover, since $\mb{\tilde{B}}(s)$ is close to being a periodic function, we find that we can exploit the regularity of its oscillations to obtain simple expressions for $\mb{\tilde{B}}_\text{reg}(s)$, and hence $\tilde{E}(s)$ and $\bar{E}_\text{tot}$. As we will see, $\mb{\tilde{B}}(s)$ can be expressed as the sum of a highly oscillatory function with slowly varying amplitude, and a smooth function that varies slowly throughout the domain of interest. Thus, the asymptotic approximation that we obtain for $\mb{\tilde{B}}(s)$ bears a strong resemblance to a 
multiple scales approximation. 

The main difficulty that we encounter in using this approach is that $\mb{\tilde{B}}(s)$ has singularities at $s = i \, n^{-1}$, where $i = 0,\,1,\,2,\,\ldots,\,n$. Because the summand in \eqref{BNondimLong} is discontinuous, we cannot immediately apply the Euler--Maclaurin summation formula to approximate the sum by an integral. Instead, we need to manipulate \eqref{BNondimLong} in order to separate the `singular terms' from the terms that can be simplified using Euler--Maclaurin summation. The `singular terms' can in turn be simplified by approximating their sum by a periodic function with regularly spaced singularities and a slowly-varying amplitude. 

This approach of separating out the singularities and then dealing with `singular terms' and `nonsingular terms' separately could potentially be applied to more general problems involving the approximation of functions defined as sums. In the problem of a chain of magnets, separating the `singular terms' (which lead to a periodic function with slowly-varying amplitude) from the `nonsingular terms' (which lead to an integral via Euler--Maclaurin summation) yields a revealing separation of `local' and `nonlocal' energy contributions, which we describe in Section \ref{S:ChainShapeEnergy}. In contrast, an analysis of the sum in \eqref{Etot_Temp} using the methods described in \cite{Hall2010} would struggle to capture both of these components in the continuum interaction energy.

Lastly, we should note a couple of assumptions implicit in the nondimensionalisation \eqref{NondimMain}. As we will see, the process of separating `singular terms' from `nonsingular terms' depends on the fact that $\mb{\tilde{r}}(s)$ and $\mb{\tilde{m}}(s)$ can be expanded in Taylor series to yield asymptotic approximations of the magnetic field due to the magnets in the neighbourhood of $\mb{\tilde{r}}(s)$. However, there are a couple of  ways in which these Taylor series could fail to yield good approximations of the effects of all nearby magnets. Firstly, and most simply, the derivatives of $\mb{\tilde{r}}(s)$ or $\mb{\tilde{m}}(s)$ could be larger than order $n^{-1}$, meaning that Taylor expansions about $s$ for $\mb{\tilde{r}}(i\,n^{-1})$ and $\mb{\tilde{m}}(i\,n^{-1})$  would no longer be asymptotic as $n \rightarrow \infty$, even when $i = s \, n + \Ord{1}$. Thus, an alternative method is required if the dipole orientation changes significantly from magnet to magnet, or if the chain is very tightly 
curved. 

Secondly, the chain could come close to itself, so that there are cases where $\mb{\tilde{r}}(j\,n^{-1}) - \mb{\tilde{r}}(s)$ is $\tOrd{n^{-1}}$, even when $j - s \, n$ is large. In this case, the magnet at $\mb{\tilde{r}}(j\,n^{-1})$ would have a significant effect on the field at $\mb{\tilde{r}}(s)$, but this would not be captured in the Taylor expansion. While a very simple modification of the method described in Section \ref{S:Main} enables us to deal with rings of magnets (as discussed in Section \ref{S:Rings}), this limitation means that we cannot approximate the energy of lattice using the one-dimensional technique presented here.

It may be noted that this restriction on the size of $\mb{\tilde{r}}(\eta) - \mb{\tilde{r}}(s)$ can be combined with the restriction on the derivatives of $\mb{\tilde{r}}(s)$ by appealing to the concept of global curvature introduced in \cite{Gonzalez1999}. In short, we find that our asymptotic approximation is valid on the condition that $||\mb{\tilde{m}}'(s)|| \ll n$ and $\rho_G \ll n$, where $\rho_G$ is the global curvature of the chain.

\section{Discrete-to-continuum analysis of a chain of spherical magnets}
\label{S:Main}

\subsection{Separation of singular and nonsingular terms in the magnetic field}
\label{S:Field}

Our first step in developing a continuum model of a chain of spherical magnets is to use discrete-to-continuum asymptotics to construct a simplified expression for $\mb{\tilde{B}}(s)$ in the case where $\mb{\tilde{r}}(s)$ and $\mb{\tilde{m}}(s)$ are specified smooth functions. From \eqref{BNondimLong}, we note that
\begin{equation}
 \mb{\tilde{B}}(s) = \sum_{i=0}^n \mb{\tilde{B}}_i(s), \label{BSimpleSum}
\end{equation} 
where
\begin{multline}
 \mb{\tilde{B}}_i(s) 
 = n^{-3}  \frac{3 \, \big\{[\mb{\tilde{r}}(s) - \mb{\tilde{r}}(i\,n^{-1})] \cdot \mb{\tilde{m}}(i\,n^{-1})\big\} \, \big[\mb{\tilde{r}}(s) - \mb{\tilde{r}}(i\,n^{-1})\big]}{||\mb{\tilde{r}}(s) - \mb{\tilde{r}}(i\,n^{-1})||^5}  \\
 - n^{-3} \frac{\mb{\tilde{m}}(i\,n^{-1})}{||\mb{\tilde{r}}(s) - \mb{\tilde{r}}(i\,n^{-1})||^3}. \label{BiNondim}
\end{multline} 
 
As we will see, the highly singular terms in \eqref{BiNondim} mean that the energy associated with `near neighbour' interactions dominates the energy due to nonlocal interactions at leading order. However, we will need to go to a higher order to be able to explore the energy effects of deforming a chain of magnets; at this higher order, both local and nonlocal terms become important. This means that the leading-order asymptotic techniques developed in \cite{Hall2010} cannot be applied to the present problem. Instead, we proceed by expanding the functions $\mb{B}_i(s)$ as Taylor series and separating the singular parts of the expressions for $\mb{B}_i(s)$ from the nonsingular parts before we take the sum given in \eqref{BSimpleSum}. This enables us to express $\mb{\tilde{B}}(s)$ as the sum of a rapidly varying and highly singular oscillatory function (representative of the local terms) and a slowly varying singular integral term (
representative of the nonlocal terms). 

The first step in our analysis is to find an asymptotic expansion for $\mb{\tilde{B}}_i(s)$ in the neighbourhood of $s = i\,n^{-1}$, concentrating especially on the terms that blow up as $s \rightarrow i \, n^{-1}$, since these are the terms that will give rise to the rapidly oscillating part of $\mb{\tilde{B}}(s)$. Towards this aim, we note that
\begin{multline}
 \mb{\tilde{r}}(s) - \mb{\tilde{r}}(i\,n^{-1}) = \mb{\tilde{r}}'(s) \, (s - i \, n^{-1}) - \frac{\mb{\tilde{r}}''(s)}{2} \, (s - i \, n^{-1})^2 \\ + \frac{\mb{\tilde{r}}'''(s) }{6} \, (s - i \, n^{-1})^3 + \Ord{(s-i\,n^{-1})^4}, \label{RDiff}
\end{multline}
and hence
\begin{multline}
 ||\mb{\tilde{r}}(s) - \mb{\tilde{r}}(i\,n^{-1})||^2 
 = ||\mb{\tilde{r}}'(s)||^2 \, (s - i \, n^{-1})^2 
 - \mb{\tilde{r}}'(s) \cdot \mb{\tilde{r}}''(s) \, (s - i \, n^{-1})^3 \\
 + \left(\frac{||\mb{\tilde{r}}''(s)||^2}{4}  + \frac{\mb{\tilde{r}}'(s) \cdot \mb{\tilde{r}}'''(s)}{3}  \right) \,  (s - i \, n^{-1})^4 + \Ord{(s - i \, n^{-1})^5 }. \label{RMagSq_FirstApprox}
\end{multline}

Moreover, substituting $i = j \pm 1$ and $i = j$ into \eqref{RDiff} and taking the difference, we find that 
\begin{multline}
 \mb{\tilde{r}}([j \pm 1]\,n^{-1}) - \mb{\tilde{r}}(j\,n^{-1}) = \pm n^{-1} \, \mb{\tilde{r}}'(s) + \frac{n^{-2}}{2} \, \mb{\tilde{r}}''(s) \pm \frac{n^{-3}}{6} \, \mb{\tilde{r}}'''(s) \\ + \Ord{n^{-4},\, (s-i\,n^{-1})^4}.
\end{multline}
Since 
$||\mb{\tilde{r}}((j\pm1)\,n^{-1} ) - \mb{\tilde{r}}(j\,n^{-1})|| = n^{-1}$, it follows that 
\begin{equation}
 1 = ||\mb{\tilde{r}}'(s)||^2 \mp  \mb{\tilde{r}}'(s) \cdot \mb{\tilde{r}}''(s) \, n^{-1} +  \left(\frac{||\mb{\tilde{r}}''(s)||^2}{4}  + \frac{\mb{\tilde{r}}'(s) \cdot \mb{\tilde{r}}'''(s)}{3}  \right) \, n^{-2} + \Ord{n^{-3}}. \label{GeometricCorrTemp}
\end{equation} 
Adding the two equations implicit in \eqref{GeometricCorrTemp} then gives
\begin{equation}
 ||\mb{\tilde{r}}'(s)||^2 = 1 - n^{-2} \, \left(\frac{||\mb{\tilde{r}}''(s)||^2}{4}  + \frac{\mb{\tilde{r}}'(s) \cdot \mb{\tilde{r}}'''(s)}{3}  \right) + \Ord{n^{-4}}, \label{rMagSq_Temp}
\end{equation}
while subtracting the two equations yields
\begin{equation}
 \mb{\tilde{r}}'(s) \cdot \mb{\tilde{r}}''(s) = \Ord{n^{-2}}. \label{r'r''ortho}
\end{equation}
Differentiating \eqref{r'r''ortho} with respect to $s$ and rearranging, we find that
\begin{equation} 
 \mb{\tilde{r}}'(s) \cdot \mb{\tilde{r}}'''(s) 
 = - ||\mb{\tilde{r}}''(s)||^2 
 + \Ord{n^{-2}},
 \label{r'''Simplification}
\end{equation}
and hence \eqref{rMagSq_Temp} yields
\begin{equation}
 ||\mb{\tilde{r}}'(s)|| = 1 + n^{-2} \, \frac{||\mb{\tilde{r}}''(s)||^2}{24} + \Ord{n^{-4}}. \label{rMag}
\end{equation}

Substituting \eqref{r'r''ortho} and \eqref{rMag} into \eqref{RMagSq_FirstApprox} and taking the square root, we obtain the following series expansion for $||\mb{\tilde{r}}(s) - \mb{\tilde{r}}(i\,n^{-1})||$:
\begin{multline}
 ||\mb{\tilde{r}}(s) - \mb{\tilde{r}}(i\,n^{-1})|| 
 = |s - i \, n^{-1}|  
 - \frac{||\mb{\tilde{r}}''(s)||^2}{24} \,  \left[ |s - i \, n^{-1}|^3 - n^{-2} \, |s - i \, n^{-1}|\right] \\
 + \Ord{(s - i \, n^{-1})^4, \, n^{-2} \, (s - i \, n^{-1})^2, \, n^{-4} \, (s-i\,n^{-1})}.
\label{RMag} 
\end{multline} 

We can now use \eqref{RDiff} and \eqref{RMag} to obtain an asymptotic expansion for $\mb{\tilde{B}}_i(s)$. Firstly, we note that
\begin{multline}
 \big[(\mb{\tilde{r}}(s) - \mb{\tilde{r}}(i\,n^{-1})) \cdot \mb{\tilde{m}}(i\,n^{-1})\big] \, (\mb{\tilde{r}}(s) - \mb{\tilde{r}}(i\,n^{-1})) 
 = (s - i\,n^{-1})^2 \, (\mb{\tilde{r}}'(s) \cdot \mb{\tilde{m}}(s)) \, \mb{\tilde{r}}'(s) \\
 - (s - i\,n^{-1})^3 \, 
 \Bigg[ \frac{\mb{\tilde{r}}''(s) \cdot \mb{\tilde{m}}(s)}{2} \, \mb{\tilde{r}}'(s) 
 + (\mb{\tilde{r}}'(s) \cdot \mb{\tilde{m}}'(s) ) \, \mb{\tilde{r}}'(s) 
 + \frac{\mb{\tilde{r}}'(s) \cdot \mb{\tilde{m}}(s)}{2} \, \mb{\tilde{r}}''(s) \Bigg] \\
 + (s - i\,n^{-1})^4 \,
 \Bigg[\frac{\mb{\tilde{r}}'''(s) \cdot \mb{\tilde{m}}(s)}{6} \, \mb{\tilde{r}}'(s)
 + \frac{\mb{\tilde{r}}'(s) \cdot \mb{\tilde{m}}''(s)}{2} \, \mb{\tilde{r}}'(s)
 + \frac{\mb{\tilde{r}}'(s) \cdot \mb{\tilde{m}}(s)}{6} \, \mb{\tilde{r}}'''(s) \\
 + \frac{\mb{\tilde{r}}''(s) \cdot \mb{\tilde{m}}'(s)}{2} \, \mb{\tilde{r}}'(s)
 + \frac{\mb{\tilde{r}}''(s) \cdot \mb{\tilde{m}}(s)}{4} \, \mb{\tilde{r}}''(s)
 + \frac{\mb{\tilde{r}}'(s) \cdot \mb{\tilde{m}}'(s)}{2} \, \mb{\tilde{r}}''(s)
 \Bigg],
\end{multline}
and that
\begin{multline}
 ||\mb{\tilde{r}}(s) - \mb{\tilde{r}}(i\,n^{-1})||^{-a} = |s - i \, n^{-1}|^{-a} 
 + a \, \frac{||\mb{\tilde{r}}''(s)||^2}{24} \,  \big( |s - i \, n^{-1}|^{-a+2} - n^{-2} \, |s - i \, n^{-1}|^{-a}\big) \\
 + \Ord{|s - i \, n^{-1}|^{-a+3}, \, n^{-2} \, |s - i \, n^{-1}|^{-a+1}, \, n^{-4} \, |s-i\,n^{-1}|^{-a}}.
\end{multline}

Substituting into \eqref{BiNondim}, it follows that $\mb{\tilde{B}}_i(s)$ has an expansion in powers of $(s-i\,n^{-1})$ of the form
\begin{equation}
 \mb{\tilde{B}}_i(s) = n^{-3} \, \frac{\mb{\Phi}_3(s)}{|s-i\,n^{-1}|^3} + n^{-3} \, \frac{\signum(s - i\,n^{-1}) \, \mb{\Phi}_2(s)}{|s-i\,n^{-1}|^2} + n^{-3} \,\frac{\mb{\Phi}_1(s)}{|s-i\,n^{-1}|} + \Ord{1}, \label{BiExpansionGeneral}
\end{equation}
where $\mb{{\Phi}}_3(s)$, $\mb{{\Phi}}_2(s)$ and $\mb{{\Phi}}_1(s)$ are given by
\begin{multline}
 \mb{\Phi}_3(s) =  3 \, (\mb{\tilde{r}}'(s) \cdot \mb{\tilde{m}}(s)) \, \mb{\tilde{r}}'(s) 
 - \mb{\tilde{m}}(s) \\
 + n^{-2} \, \Bigg(
 - \frac{5 \, ||\mb{\tilde{r}}''(s)||^2 \, (\mb{\tilde{r}}'(s) \cdot \mb{\tilde{m}}(s))}{8} \, \mb{\tilde{r}}'(s)  
 + \frac{||\mb{\tilde{r}}''(s)||^2 }{8}\, \mb{\tilde{m}}(s) 
 \Bigg) 
 + \Ord{n^{-4}}, \label{Phi3Defn}
\end{multline}
\begin{multline}
 \mb{\Phi}_2(s) = 
 - \frac{3 \, \mb{\tilde{r}}''(s) \cdot \mb{\tilde{m}}(s)}{2} \, \mb{\tilde{r}}'(s) 
 - 3 \, (\mb{\tilde{r}}'(s) \cdot \mb{\tilde{m}}'(s) ) \, \mb{\tilde{r}}'(s)  \\
 - \frac{3\, \mb{\tilde{r}}'(s) \cdot \mb{\tilde{m}}(s)}{2} \, \mb{\tilde{r}}''(s) 
 + \mb{\tilde{m}}'(s) + \Ord{n^{-2}},
\end{multline}
and
\begin{multline}
 \mb{\Phi}_1(s) =
 \frac{\mb{\tilde{r}}'''(s) \cdot \mb{\tilde{m}}(s)}{2} \, \mb{\tilde{r}}'(s)
 + \frac{3 \, \mb{\tilde{r}}'(s) \cdot \mb{\tilde{m}}''(s)}{2} \, \mb{\tilde{r}}'(s)
 + \frac{\mb{\tilde{r}}'(s) \cdot \mb{\tilde{m}}(s)}{2} \, \mb{\tilde{r}}'''(s) \\
 + \frac{3 \,\mb{\tilde{r}}''(s) \cdot \mb{\tilde{m}}'(s)}{2} \, \mb{\tilde{r}}'(s)
 + \frac{3 \, \mb{\tilde{r}}''(s) \cdot \mb{\tilde{m}}(s)}{4} \, \mb{\tilde{r}}''(s)
 + \frac{3 \, \mb{\tilde{r}}'(s) \cdot \mb{\tilde{m}}'(s)}{2} \, \mb{\tilde{r}}''(s) \\
 + \frac{5 \, ||\mb{\tilde{r}}''(s)||^2 \, (\mb{\tilde{r}}'(s) \cdot \mb{\tilde{m}}(s))}{8} \, \mb{\tilde{r}}'(s) 
 - \frac{1}{2} \, \mb{\tilde{m}}''(s)
 - \frac{||\mb{\tilde{r}}''(s)||^2}{8} \, \mb{\tilde{m}}(s)
 + \Ord{n^{-2}}. \label{Phi1Defn}
\end{multline}

When we take a sum of \eqref{BiExpansionGeneral} over all $i$, we will find that $\mb{\Phi}_3(s)$, $\mb{\Phi}_2(s)$ and $\mb{\Phi}_1(s)$ are the slowly-varying amplitudes of some highly singular oscillating functions. In contrast, the order one terms in \eqref{BiExpansionGeneral} are nonsingular, and we can therefore replace the sum of these order one terms with an integral. With this in mind, we rearrange \eqref{BSimpleSum} by adding and subtracting sums of singular terms to obtain
\begin{equation}
 \mb{\tilde{B}}(s) = \mb{S}_*(s)+ \sum_{k=1}^3 \mb{\Phi}_k(s) \, S_k(s), 
 \label{BInTermsOfS}
\end{equation} 
where
\begin{equation}
 S_k(s) = \sum_{i=0}^n \frac{\signum(s - i \, n^{-1}) \, n^{-3}}{(s - i \, n^{-1})^{k}}, \label{SkDefn}
\end{equation}
and
\begin{equation}
 \mb{S}_*(s) = \sum_{i=0}^n \left(\mb{\tilde{B}}_i(s) 
 - \frac{n^{-3} \, \mb{\Phi}_3(s)}{|s-i\,n^{-1}|^3} 
 - \frac{\signum(s - i\,n^{-1}) \, n^{-3} \, \mb{\Phi}_2(s)}{|s-i\,n^{-1}|^2} 
 - \frac{n^{-3} \, \mb{\Phi}_1(s)}{|s-i\,n^{-1}|}\right).
\end{equation}


\subsection{Approximation of sums}
\label{S:ApproxSums}

Having separated the slowly-varying behaviour of $\mb{\tilde{B}}(s)$ from the highly singular oscillatory behaviour, we can now approximate the sums in \eqref{BInTermsOfS} using classical asymptotic methods. Importantly, the summand of $\mb{S}_*(s)$ is continuous, whereas the summand in our original definition of $\mb{\tilde{B}}(s)$ was not. Thus, we can now use the Euler--Maclaurin summation formula (see, for example, \cite{CKP}) to approximate $\mb{S}_*(s)$ with an integral as follows:
\begin{multline}
 \mb{S}_*(s) = n^{-2} \, \int_0^1 \Bigg[ \frac{3 \, \big[(\mb{\tilde{r}}(s) - \mb{\tilde{r}}(\eta)) \cdot \mb{\tilde{m}}(\eta)\big] \, (\mb{\tilde{r}}(s) - \mb{\tilde{r}}(\eta))}{||\mb{\tilde{r}}(s) - \mb{\tilde{r}}(\eta)||^5}  
 - \frac{\mb{\tilde{m}}(\eta)}{||\mb{\tilde{r}}(s) - \mb{\tilde{r}}(\eta)||^3} \\
 - \frac{\mb{\Phi}_3(s)}{|s-\eta|^3} 
 - \frac{\signum(s - \eta) \, \mb{\Phi}_2(s)}{|s-\eta|^2} 
 - \frac{\mb{\Phi}_1(s)}{|s-\eta|}
 \Bigg] \, d\eta
 + \Ord{n^{-3}}. \label{S*Simp}
\end{multline} 

Having approximated $\mb{S}_*(s)$ with an integral, we now wish to approximate $S_1(s)$, $S_2(s)$ and $S_3(s)$, defined in \eqref{SkDefn}, by noting that they are close to being periodic functions in a fast variable, $X = s\,n$. One easy way to achieve this multiple scales approximation of the functions $S_k(s)$ is to replace the finite sums in \eqref{SkDefn} with infinite sums. For $S_3$, for example, we find that
\begin{equation}
 S_3(s,\,X) = \sum_{i=-\infty}^{\infty} \frac{1}{|X-i|^3} - n^{-3} \, \sum_{i=-\infty}^{-1} \frac{1}{|s - i \, n^{-1}|^{3}} - n^{-3} \, \sum_{i=n+1}^{\infty} \frac{1}{|s - i \, n^{-1}|^{-3}}. \label{S3SimplifiedSums}
\end{equation}
Since $0 < s < 1$, the latter two sums in \eqref{S3SimplifiedSums} can be approximated by integrals using Euler--Maclaurin series as long as $s \gg n^{-1}$ and $1 - s \gg n^{-1}$. This yields
\begin{equation}
  S_3(s,\,X) = \Lambda_3(X) - n^{-2}\, \int_{-\infty}^0 \frac{1}{|s - \eta|^3} \, d\eta - n^{-2} \, \int_{1}^\infty \frac{1}{|s - \eta|^3} \, d\eta + \Ord{n^{-3}},  \label{S3SimplifiedTemp}
\end{equation} 
where $\Lambda_3(X)$ is the periodic function defined by the convergent sum
\begin{equation}
 \Lambda_3(X) = \sum_{i=-\infty}^{\infty} \frac{1}{|X-i|^3}. \label{Lambda3Def}
\end{equation} 
Care needs to be taken with $S_3$ (and similarly $S_2$ and $S_1$) if $s = \Ord{n^{-1}}$ or if $1 - s = \Ord{n^{-1}}$. This is because the singularity in the summand causes conventional Euler--Maclaurin summation to yield series that are not asymptotic for the latter sums in \eqref{S3SimplifiedSums}. This problem indicates that we need to consider a discrete problem in boundary layer regions at the ends of the chain, analogous to the boundary layers described in \cite{Hall2011,Hall2010}.

The two integrals in \eqref{S3SimplifiedTemp} can both be evaluated analytically, giving the result that
\begin{align}
 -\int_{-\infty}^0 \frac{1}{|s - \eta|^3} \, d\eta - \int_{1}^\infty \frac{1}{|s - \eta|^3} \, d\eta 
 &= -\frac{1}{2 \, s^{2}} - \frac{1}{2 \, (1- s)^2} \notag \\
 &= \dashint_0^1 \frac{1}{|s - \eta|^3} \, d\eta, \label{FirstFinitePartInt}
\end{align}
where the dashed integral sign is used to indicate a finite part integral (see, for example, \cite{RichardsDistributions}). 


Using this result, we find that \eqref{S3SimplifiedTemp} becomes
\begin{equation}
 S_3(s,\,X) = \Lambda_3(X) + n^{-2} \, \dashint_0^1 \frac{1}{|s - \eta|^3} \, d\eta + \Ord{n^{-3}}. \label{S3Simp}
\end{equation} 
Similarly, it can be shown that
\begin{equation}
 S_2(s,\,X) = n^{-1} \, \Lambda_2(X) + n^{-2} \, \dashint_{0}^1 \frac{\signum(s - \eta)}{|s - \eta|^2} \, d\eta + \Ord{n^{-3}}, \label{S2Simp}
\end{equation} 
where
\begin{equation}
 \Lambda_2(X) = \sum_{-\infty}^\infty \frac{\signum(X - i)}{|X - i|^2}. \label{Lambda2Def}
\end{equation}

The infinite sum associated with $S_1$ requires a bit more care, however, because of the divergence of the harmonic series. We note that $S_1$ can be rewritten in the form
\begin{multline}
 S_1(s,\,X) = \underbrace{n^{-2}\, \sum_{i=-K}^{K} \frac{1}{|X-i|} 
 - 2 \, n^{-2} \, \log K}_{S_{1\text{a}}} 
 + 2 \, n^{-2} \, \log n \\
 - n^{-2} \, \Bigg( \underbrace{n^{-1}\,\sum_{i=-K}^{-1} \frac{1}{|s - i \, n^{-1}|} - \log \left(K \, n^{-1}\right)}_{S_{1\text{b}}} \Bigg) \\
 - n^{-2} \, \Bigg( \underbrace{n^{-1}\,\sum_{i=n+1}^{K} \frac{1}{|s - i \, n^{-1}|} - \log \left(K \, n^{-1}\right)}_{S_{1\text{c}}} \Bigg). 
\end{multline}
This representation is chosen so that $S_{1\text{a}}$ is convergent in the limit $K \rightarrow \infty$ and so that $S_{1\text{b}}$ and $S_{1\text{c}}$ can be used to construct a finite part integral as before. Proceeding as above, we find that
\begin{equation}
 S_1(s,\,X) = 2 \, n^{-2} \, \log n + n^{-2} \, \Lambda_1(X) + n^{-2} \dashint_0^1 \frac{1}{|s - \eta|} \, d\eta + \Ord{n^{-3}}, \label{S1Simp}
\end{equation} 
where
\begin{equation}
 \Lambda_1(X) = \lim_{K \rightarrow \infty} \left[-2 \, \log K +  \sum_{i=-K}^K \frac{1}{|X - i|}\right]. \label{Lambda1Def}
\end{equation}

Now we can substitute \eqref{S*Simp}, \eqref{S3Simp}, \eqref{S2Simp} and \eqref{S1Simp} into \eqref{BInTermsOfS}, noting that the finite part integrals in \eqref{S3Simp}, \eqref{S2Simp} and \eqref{S1Simp} will cancel with parts of the integral in \eqref{S*Simp}. This yields the following expression for $\mb{\tilde{B}}(s)$:
\begin{multline}
  \mb{\tilde{B}}(s) = \mb{\Phi}_3(s) \, \Lambda_3(X) 
  + n^{-1} \, \mb{\Phi}_2(s) \, \Lambda_2(X) 
  + 2 \, n^{-2} \, \log n \, \mb{\Phi}_1(s) 
  + n^{-2} \, \mb{\Phi}_1(s) \, \Lambda_1(X) \\
  +  n^{-2} \, \dashint_0^1 \frac{3 \, \big[(\mb{\tilde{r}}(s) - \mb{\tilde{r}}(\eta)) \cdot \mb{\tilde{m}}(\eta)\big] \, (\mb{\tilde{r}}(s) - \mb{\tilde{r}}(\eta))}{||\mb{\tilde{r}}(s) - \mb{\tilde{r}}(\eta)||^5}  
 - \frac{\mb{\tilde{m}}(\eta)}{||\mb{\tilde{r}}(s) - \mb{\tilde{r}}(\eta)||^3 } \, d\eta + \Ord{n^{-3}}. \label{BFinal}
\end{multline}

We now have an equation for $\mb{\tilde{B}}(s)$  that does not involve any sums and that depends only on the continuum functions $\mb{\tilde{r}}(s)$ and $\mb{\tilde{m}}(s)$. Importantly, we find that the rapidly varying and almost periodic behaviour of the magnetic field -- represented by the $\Lambda_i(X)$ terms -- is separated from the long-scale changes in the magnetic field represented by the singular integral. 

\subsection{Energy of a chain of magnets}

In order to calculate the energy defined in \eqref{EDefn}, it is useful to have a simple expression for $\mb{\tilde{B}}_\text{reg}(s)$. To this end, we note that substituting \eqref{BiExpansionGeneral} and \eqref{BFinal} into \eqref{BRegDefn} yields
\begin{multline}
 \mb{\tilde{B}}_\text{reg}(i \,n^{-1})
 = \lim_{s\rightarrow i\,n^{-1}} 
 \Bigg[ 
 \mb{\Phi}_3(s) \, \left(\Lambda_3(s \, n) - \frac{1}{|s \, n - i|^3}\right) \\
 + n^{-1} \, \mb{\Phi}_2(s) \, \left(\Lambda_2(s \, n) - \frac{\signum(s \, n - i)}{|s \, n - i|^2}\right) 
 + 2 \, n^{-2} \, \log n \, \mb{\Phi}_1(s) \\
 + n^{-2} \, \mb{\Phi}_1(s) \, \left(\Lambda_1(s \, n) - \frac{1}{|s \, n - i|}\right) \\
 +  n^{-2} \, \dashint_0^1 \frac{3 \, \big[(\mb{\tilde{r}}(s) - \mb{\tilde{r}}(\eta)) \cdot \mb{\tilde{m}}(\eta)\big] \, (\mb{\tilde{r}}(s) - \mb{\tilde{r}}(\eta))}{||\mb{\tilde{r}}(s) - \mb{\tilde{r}}(\eta)||^5}  
 - \frac{\mb{\tilde{m}}(\eta)}{||\mb{\tilde{r}}(s) - \mb{\tilde{r}}(\eta)||^3 } \, d\eta
 \Bigg]
 + \Ord{n^{-3}}.
\end{multline}
From the definitions of  $\Lambda_3(X)$, $\Lambda_2(X)$ and $\Lambda_1(X)$ in \eqref{Lambda3Def}, \eqref{Lambda2Def} and  \eqref{Lambda1Def}, it follows that
\begin{gather}
 \lim_{X \rightarrow i} \left(\Lambda_3(X) - \frac{1}{|X - i|^3}\right) = \sum_{\substack{i=-\infty \\ i \neq 0}}^{\infty} \frac{1}{|i|^{3}} = 2 \, \zeta(3),
\\ 
 \lim_{X \rightarrow i} \left(\Lambda_2(X) - \frac{\signum(X - i)}{|X - i|^2}\right) = \sum_{\substack{i=-\infty \\ i \neq 0}}^{\infty} \frac{\signum(i)}{i^{2}} = 0,
\\ 
 \lim_{X \rightarrow i} \left(\Lambda_1(X) - \frac{1}{|X - i|}\right) = \lim_{K \rightarrow \infty}\left[\sum_{\substack{i=-K \\ i \neq 0}}^{K} \frac{1}{|i|} - 2 \, \log K\right] = 2 \, \gamma,
\end{gather}
where $\zeta(z)$ is the Riemann Zeta function and $\gamma$ is the Euler--Mascheroni constant.
Hence, it is convenient for us to define $\mb{\tilde{B}}_\text{reg}(s)$ so that
\begin{multline}
 \mb{\tilde{B}}_\text{reg}(s) 
 = 
 2 \, \zeta(3) \, \mb{\Phi}_3(s)
 + 2 \, \mb{\Phi}_1(s) \, n^{-2} \, \log n
 + 2 \, \gamma \, \mb{\Phi}_1(s) \, n^{-2} \\
 +  n^{-2} \, \dashint_0^1 \frac{3 \, \big[(\mb{\tilde{r}}(s) - \mb{\tilde{r}}(\eta)) \cdot \mb{\tilde{m}}(\eta)\big] \, (\mb{\tilde{r}}(s) - \mb{\tilde{r}}(\eta))}{||\mb{\tilde{r}}(s) - \mb{\tilde{r}}(\eta)||^5}  
 - \frac{\mb{\tilde{m}}(\eta)}{||\mb{\tilde{r}}(s) - \mb{\tilde{r}}(\eta)||^3 } \, d\eta \label{BregAsymp}
 + \Ord{n^{-3}},
\end{multline}

Given this expression for $\mb{\tilde{B}}_\text{reg}(s)$, we can easily obtain asymptotic expressions for the energy of any individual magnet and for the entire chain. Substituting \eqref{BregAsymp} into \eqref{EDefn} and exploiting \eqref{mNorm=1}, \eqref{Phi3Defn} and \eqref{Phi1Defn}, we find that
\begin{multline}
 \tilde{E}(s) = - 2 \, \zeta(3) \, \big\{ 3 \, [\mb{\tilde{r}}'(s) \cdot \mb{\tilde{m}}(s)]^2 - 1\big\} \\
 + n^{-2} \, \log n \, \Bigg[
 - 2 \, (\mb{\tilde{r}}'''(s) \cdot \mb{\tilde{m}}(s)) \, (\mb{\tilde{r}}'(s) \cdot \mb{\tilde{m}}(s))
 - 3 \, (\mb{\tilde{r}}'(s) \cdot \mb{\tilde{m}}''(s)) \, (\mb{\tilde{r}}'(s) \cdot \mb{\tilde{m}}(s)) \\
 - 3 \, (\mb{\tilde{r}}''(s) \cdot \mb{\tilde{m}}'(s))\, (\mb{\tilde{r}}'(s) \cdot \mb{\tilde{m}}(s)) 
 - \frac{3 \, (\mb{\tilde{r}}''(s) \cdot \mb{\tilde{m}}(s))^2}{2}
 - 3 \, (\mb{\tilde{r}}'(s) \cdot \mb{\tilde{m}}'(s)) \, (\mb{\tilde{r}}''(s) \cdot \mb{\tilde{m}}(s))  \\
 - \frac{5 \, ||\mb{\tilde{r}}''(s)||^2 \, (\mb{\tilde{r}}'(s) \cdot \mb{\tilde{m}}(s))^2}{4}
 + \mb{\tilde{m}}''(s) \cdot \mb{\tilde{m}}(s)
 + \frac{||\mb{\tilde{r}}''(s)||^2}{4}
 \Bigg]
 + \Ord{n^{-2}}. \label{EnergyAsympOrig}
\end{multline} 

\subsection{Dipoles align tangential to the chain}

In the absence of an applied magnetic field, we intuitively expect the dipoles of the spherical magnets to line up with the tangent of the chain. However, it is reassuring to note that we can also obtain this result from first principles. Our approach is to assume that $\mb{\tilde{r}}(s)$ is specified and to seek a solution for $\mb{\tilde{m}}(s)$ that minimizes the total energy.
This will ultimately give us $\mb{\tilde{m}}(s)$ as an asymptotic series in $n$ that depends on $\mb{\tilde{r}}(s)$. We can then substitute this solution for $\mb{\tilde{m}}(s)$ into \eqref{EnergyAsympOrig} in order to obtain a formula for the total energy as a function of $\mb{\tilde{r}}(s)$ alone.

Given the form of \eqref{EnergyAsympOrig}, it seems appropriate to expand $\mb{\tilde{m}}(s)$ as an asymptotic series as follows:
\begin{equation}
 \mb{\tilde{m}}(s) = \mb{\tilde{m}}_0(s) + n^{-2} \, \log n \, \mb{\tilde{m}}_1(s) + n^{-2} \, \mb{\tilde{m}}_2(s) + \Ord{n^{-3}}. \label{MAsymp}
\end{equation}
Since $||\mb{\tilde{m}}(s)|| = 1$, we further note that
\begin{equation}
 1 = ||\mb{\tilde{m}}_0(s)||^2 + 2 \, n^{-2} \, \log n \, \mb{\tilde{m}}_0(s) \cdot \mb{\tilde{m}}_1(s) + 2 \, n^{-2} \, \mb{\tilde{m}}_0(s) \cdot \mb{\tilde{m}}_2(s) + \Ord{n^{-3}}. \label{MMagAsymp}
\end{equation}
Without loss of generality, we choose $\mb{\tilde{m}}_0(s)$ so that $||\mb{\tilde{m}}_0(s)||^2 = 1$, we choose $\mb{\tilde{m}}_1(s)$ so that $\mb{\tilde{m}}_0(s) \cdot \mb{\tilde{m}}_1(s) = 0$, and we choose  $\mb{\tilde{m}}_2(s)$ so that $\mb{\tilde{m}}_0(s) \cdot \mb{\tilde{m}}_2(s) = 0$.

Substituting into \eqref{EnergyAsympOrig}, this yields the result that
\begin{equation}
 \tilde{E}(s) = \underbrace{- 2 \, \zeta(3) \, \big[ 3 \, (\mb{\tilde{r}}'(s) \cdot \mb{\tilde{m}}_0(s))^2 - 1\big]}_{\tilde{E}_0(s)} + \Ord{n^{-2} \, \log n}. \label{ELeading}
\end{equation}
If we now pick $\mb{\tilde{m}}_0(s)$ so that the $\tilde{E}_0(s)$ is minimized for all $s$, it will follow that the total energy is minimized to leading order. Simple inspection of \eqref{ELeading} indicates that $\tilde{E}_0(s)$ is minimized when $\mb{\tilde{m}}_0(s)$ and $\mb{\tilde{r}}'(s)$ are parallel. Thus, we recover the expected result that the dipole moments are always aligned parallel to the tangent of the chain.

We note that the energy in \eqref{ELeading} would be maximized by taking $\mb{\tilde{r}}'(s) \cdot \mb{\tilde{m}}_0(s) = 0$. That is, our formulation predicts the maximum energy when the dipoles are perpendicular to the tangent of the chain. This may be surprising, since intuition suggests that the energy maximum would be achieved when the dipoles are aligned with the tangent of the chain, but where they alternate in sign. However, this (and many other possible critical points in the energy surface) are excluded from our analysis by the assumption that $\tilde{m}$ is a slowly varying function of $s$. Hence, we have only demonstrated that dipoles align with the tangent of the chain in the case where there are no rapid changes in dipole orientation along the chain.

\subsection{Relationship between chain shape and energy}
\label{S:ChainShapeEnergy}

Having demonstrated that dipoles will align themselves parallel to the tangent of the chain, we now assume, without loss of generality, that they point along the tangent in the direction of increasing $s$. Noting the constraint on the magnitude of $\mb{\tilde{m}}_0(s)$, we find that
\begin{equation}
 \mb{\tilde{m}}_0(s) = \frac{\mb{\tilde{r}}'(s)}{||\mb{\tilde{r}}'(s)||}. 
\end{equation}
Using \eqref{rMag}, it follows that 
\begin{equation}
 \mb{\tilde{m}}_0(s) = \mb{\tilde{r}}'(s) - n^{-2} \, \frac{||\mb{\tilde{r}}''(s)||^2 \, \mb{\tilde{r}}'(s)}{24} + \Ord{n^{-4}}, \label{MLeading}
\end{equation}
and we also find that
\begin{equation}
 \mb{\tilde{m}}_1(s) \cdot \mb{\tilde{r}}'(s) = \Ord{n^{-2}},
\end{equation}
and
\begin{equation}
 \mb{\tilde{m}}_2(s) \cdot \mb{\tilde{r}}'(s) = \Ord{n^{-2}}.
\end{equation}

Substituting 
\begin{equation}
 \mb{\tilde{m}}(s) = \mb{\tilde{r}}'(s) + (n^{-2} \, \log n) \, \mb{\tilde{m}}_1(s) + \Ord{n^{-2}}
\end{equation} 
into \eqref{EnergyAsympOrig} and using the fact that $\mb{\tilde{r}}'(s) \cdot \mb{\tilde{r}}''(s) = \Ord{n^{-2}}$ and $\mb{\tilde{r}}'(s) \cdot \mb{\tilde{m}}_1(s) = \Ord{n^{-2}}$, we obtain the surprising result that
\begin{equation}
 \tilde{E}(s) = - 4 \, \zeta(3) + \Ord{n^{-2}}.
\end{equation}
All of the $\Ord{n^{-2}\,\log n}$ terms have been eliminated, and we need to go to $\tOrd{n^{-2}}$ in order to obtain a term dependent on $\mb{\tilde{r}}(s)$.

Usefully, this mass cancellation occurred because $\mb{\Phi}_1(s) \cdot \mb{\tilde{m}}(s) = \tOrd{n^{-2} \, \log n}$, and it quickly follows that $\mb{\tilde{m}}_1(s) \equiv 0$, and hence $\mb{\Phi}_1(s) \cdot \mb{\tilde{m}}(s) = \tOrd{n^{-2}}$. Noting that
\begin{equation}
 \mb{\tilde{m}}(s) 
 = \mb{\tilde{r}}'(s) 
 + n^{-2} \left(
 -\frac{||\mb{\tilde{r}}''(s)||^2 \, \mb{\tilde{r}}'(s)}{24}
 + \mb{\tilde{m}}_2(s) 
 \right)
 + \Ord{n^{-3}},
\end{equation}
and that
\begin{equation}
 \mb{\tilde{m}}(s) \cdot \mb{\tilde{r}}'(s) 
 = 1 
 + n^{-2} \, \frac{||\mb{\tilde{r}}''(s)||^2}{24}
 + \Ord{n^{-3}}
\end{equation}
we can substitute into the full expansion for $\mb{\tilde{B}}_\text{reg}(i \, n^{-1}) \cdot \mb{\tilde{m}}(s)$ to show that
\begin{multline}
 \tilde{E}(s) = - 4 \, \zeta(3) 
 + n^{-2} \, \frac{\zeta(3)}{2} \, ||\mb{\tilde{r}}''(s)||^2  \\
 -  n^{-2} \, \dashint_0^1 
 \frac{3 \, \big[(\mb{\tilde{r}}(s) - \mb{\tilde{r}}(\eta)) \cdot \mb{\tilde{r}}'(\eta)\big] \, \big[(\mb{\tilde{r}}(s) - \mb{\tilde{r}}(\eta)) \cdot \mb{\tilde{r}}'(s)\big]}{||\mb{\tilde{r}}(s) - \mb{\tilde{r}}(\eta)||^5}  
 - \frac{\mb{\tilde{r}}'(\eta) \cdot \mb{\tilde{r}}'(s)}{||\mb{\tilde{r}}(s) - \mb{\tilde{r}}(\eta)||^3 } \, d\eta \\
 + \Ord{n^{-3}}. \label{EnergyAsympSimp}
\end{multline}

Ignoring the importance of boundary layers at the ends of the domain and applying Euler--Maclaurin summation, it therefore follows that
\begin{multline}
 \bar{E}_\text{tot} = - 2 \, \zeta(3) \, n - 2 \, \zeta(3) 
  + n^{-1} \, \Bigg[ \frac{\zeta(3)}{4} \int_0^{1}  \, ||\mb{\tilde{r}}''(s)||^2 \, ds \\
 - \int_0^1 \dashint_0^1 
 \frac{3 \, \big[(\mb{\tilde{r}}(s) - \mb{\tilde{r}}(\eta)) \cdot \mb{\tilde{r}}'(\eta)\big] \, \big[(\mb{\tilde{r}}(s) - \mb{\tilde{r}}(\eta)) \cdot \mb{\tilde{r}}'(s)\big]}{2 \, ||\mb{\tilde{r}}(s) - \mb{\tilde{r}}(\eta)||^5}  
 - \frac{\mb{\tilde{r}}'(\eta) \cdot \mb{\tilde{r}}'(s)}{2 \, ||\mb{\tilde{r}}(s) - \mb{\tilde{r}}(\eta)||^3 } \, d\eta \, ds \Bigg] \\
 + \Ord{n^{-2}}. \label{Etot_DtoC}
\end{multline}

In practical problems, we seek a choice of $\mb{\tilde{r}}(s)$ that minimizes $E_\text{tot}$ subject to certain constraints (or subject to an additional energy contribution due to, for example, gravity). Since the $\tOrd{n}$ term in \eqref{Etot_DtoC} is independent of $\mb{\tilde{r}}(s)$, it can be thought of as an energetic ground state, which arises because the spheres are arranged in a chain rather than separated. Instead, the relationship between $\mb{\tilde{r}}(s)$ and $E_\text{tot}$ (\ie the energetic response of the chain to deformation) is, to leading order, entirely determined by the $\tOrd{n^{-1}}$ term.

We observe that the $\tOrd{n^{-1}}$ term in \eqref{Etot_DtoC} consists of two distinct parts. Firstly, there is a local energy that is mathematically equivalent to the bending stiffness of an elastic rod (see, for example, \cite{LoveElasticity}). This reflects the fact that bending a chain of magnets causes the magnetic dipoles to be misaligned with each other, increasing the energy of the chain in direct proportion to the square of the curvature. However, there is also a nonlocal energy contribution that appears at the same asymptotic order in $n$ as this bending stiffness term. This represents the fact that each magnet feels the effective magnetic 
field of all of the other magnets in the chain, not just their nearest neighbours. The fact that both local and nonlocal interactions are important means that analysing the deformation of a chain of magnets is more complicated than analysing the deformation of an elastic rod.

An additional problem with the nonlocal term in \eqref{Etot_DtoC} is that it is not well-behaved at the ends of the chain. Consider, for example, a finite straight chain where $\mb{\tilde{r}}(s) = s \, \mb{i}$. Substituting into \eqref{EnergyAsympSimp}, we find that
\begin{equation}
 \tilde{E}(s) = -4 \, \zeta(3) - 2 \, \dashint_0^1 \frac{1}{|s - \eta|^3} \, d\eta = -4 \, \zeta(3) + n^{-2} \left(\frac{1}{s^2} + \frac{1}{ (1-s)^2} \right),
\end{equation}
and hence the integral in \eqref{Etot_DtoC} does not converge. This is a result of the fact that \eqref{EnergyAsympSimp} and \eqref{Etot_DtoC} are only appropriate when $s$ and $1 - s$ are both much larger than $n^{-1}$. As $s$ appproaches zero or one, it becomes necessary to consider the discrete boundary layer problem, and the true total energy should include contributions based on the behaviour in the boundary layers as well as a continuum energy based on \eqref{EnergyAsympSimp}. Despite this problem, it is still possible to obtain useful results from \eqref{Etot_DtoC} by restricting our attention to configurations of chains without ends. The most natural of these is a finite ring.

\section{Rings of magnets}
\label{S:Rings}

\subsection{An undeformed circular ring}

In Section \ref{S:ContinuumFormulation}, we noted that our approximation may encounter difficulties if the chain gets too close to itself. For a smooth ring, however, this is not the problem that it might appear to be: the geometry and topology of a ring mean that it is natural for us to introduce periodic extensions of $\mb{\tilde{r}}(s)$ and $\mb{\tilde{m}}(s)$ so that neighbouring points on the ring have neighbouring $s$ values.

Specifically, we consider the case where $n$ magnets are arranged smoothly in a ring so that $\mb{\tilde{r}}(0) = \mb{\tilde{r}}(1)$ represents the dimensionless location of the $n$th magnet, and  $\mb{\tilde{r}}(s)$ is a smooth periodic function of $s$ with period 1. As a result of this, the singular integral in \eqref{EnergyAsympSimp} will no longer be singular at $s=0$ and $s=1$, and hence the double integral in \eqref{Etot_DtoC} will exist.

It should be noted that \eqref{Etot_DtoC} was derived for the case of $n +1$ magnets in a chain, and we now wish to consider a ring containing only $n$ magnets. This requires only minor adjustments to our final use of the Euler--Maclaurin summation formula, and we find that the equivalent of \eqref{Etot_DtoC} for a ring of $n$ magnets is
\begin{multline}
 \bar{E}_\text{tot} = - 2 \, \zeta(3) \, n 
  + n^{-1} \, \Bigg[ \frac{\zeta(3)}{4} \int_0^{1}  \, ||\mb{\tilde{r}}''(s)||^2 \, ds \\
 - \int_0^1 \dashint_0^1 
 \frac{3 \, \big[(\mb{\tilde{r}}(s) - \mb{\tilde{r}}(\eta)) \cdot \mb{\tilde{r}}'(\eta)\big] \, \big[(\mb{\tilde{r}}(s) - \mb{\tilde{r}}(\eta)) \cdot \mb{\tilde{r}}'(s)\big]}{2 \, ||\mb{\tilde{r}}(s) - \mb{\tilde{r}}(\eta)||^5}  
 - \frac{\mb{\tilde{r}}'(\eta) \cdot \mb{\tilde{r}}'(s)}{2 \, ||\mb{\tilde{r}}(s) - \mb{\tilde{r}}(\eta)||^3 } \, d\eta \, ds \Bigg] \\
 + \Ord{n^{-2}}. \label{Etot_DtoCRing}
\end{multline}

Throughout the rest of this section, we use \eqref{Etot_DtoCRing} to approximate the energy of a perfectly circular ring of magnets, and then consider how this energy will change if the circle is slightly deformed. From this, we are able to characterize the vibrational modes of a circular ring of magnets, and hence compare them with the vibrational modes of a classical elastic ring.

The position vector for an undeformed circle of magnets can be described by
\begin{equation}
 \mb{\tilde{r}}(s) = R \, \big[\cos (2 \, \pi \, s) \, \mb{i} + \sin (2 \, \pi \, s) \, \mb{j} \big], \label{UndeformedCirc}
\end{equation}
where \eqref{NondimGapNeighbours} implies that
\begin{equation}
 R = \frac{n^{-1}}{2 \sin(\pi \, n^{-1})} = \frac{1}{2 \, \pi} + \frac{\pi}{12} \, n^{-2} + \Ord{n^{-4}}. \label{RDefn}
\end{equation}

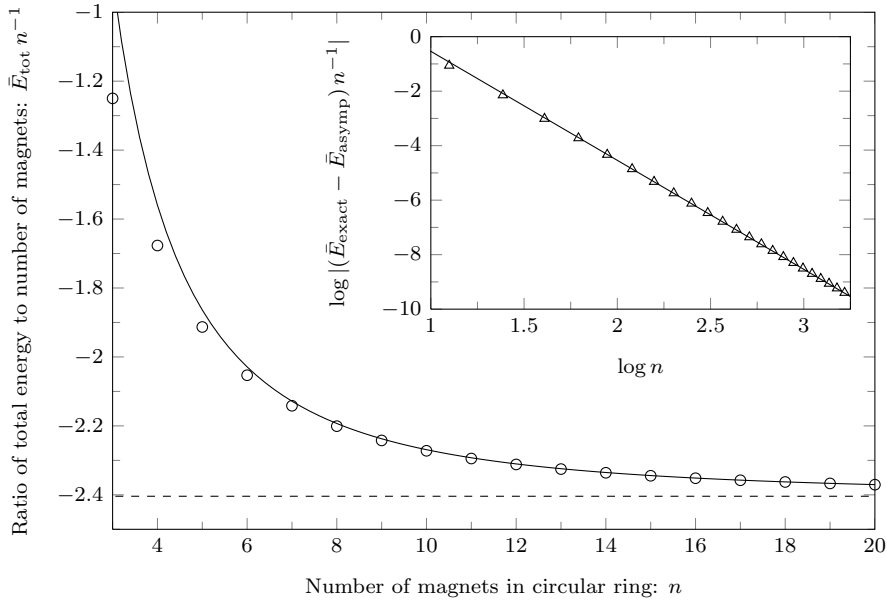
\begin{figure}[htp]
\centering
\begin{tikzpicture}
\pgfplotstableread{data-286.txt}\tableA
\pgfplotstableread{data-287.txt}\tableB
\pgfplotstableread{data-290.txt}\tableC
\pgfplotstableread{data-291.txt}\tableD
\begin{axis}[
ymin=-2.5,ymax=-1,xmin=3,xmax=20,
xlabel={Number of magnets in circular ring: $n$},
ylabel={Ratio of total energy to number of magnets: $\bar{E}_\text{tot} \, n^{-1}$},
minor tick num=1,
font=\footnotesize,
width={0.9\textwidth},
height={0.65\textwidth}
]
\addplot[color=black, only marks, mark=o] table[x index = 0, y index = 1] from \tableA;
\addplot[color=black] table[x index = 0, y index = 1] from \tableB;
\addplot[color=black, dashed] coordinates{
 (2, -2.404113806319189)
 (20, -2.404113806319189)
};
\end{axis}
\begin{axis}[
xshift=.325\textwidth,yshift=.225\textwidth,width=0.55\textwidth,height=.4\textwidth,
ymin=-10,ymax=0,xmin=1,xmax=3.25,
xlabel={$\log n$},
ylabel={$\log |(\bar{E}_\text{exact} - \bar{E}_\text{asymp}) \, n^{-1}|$},
minor tick num=1,
font=\footnotesize
]
\addplot[color=black, only marks, mark=triangle] table[x index = 0, y index = 1] from \tableC;
\addplot[color=black] table[x index = 0, y index = 1] from \tableD;
\end{axis}
\end{tikzpicture}
\caption{Main figure: Comparison of numerical calculations of the energy per magnet in a perfectly circular ring (open circles) with the asymptotic approximation obtained from \eqref{Etot_PerfectCirc} (continuous curve). As $n \rightarrow \infty$, we find that the energy per magnet tends towards the asymptotic value of $-2 \, \zeta(3)$ (dashed line). Inset figure: The discrepancy between the numerically obtained energy and the asymptotic approximation as a function of the number of magnets (triangles) shown on a logarithmic plot. The gradient of the line of best fit (continuous line) is $-4$, indicating that the next correction to \eqref{Etot_PerfectCirc} is in fact $\Ord{n^{-3}}$ rather than the $\Ord{n^{-2}}$ expected.}
\label{F:CircleEnergy}
\end{figure}

Substituting into \eqref{Etot_DtoCRing}, we find that 
\begin{multline}
 \bar{E}_\text{tot} = - 2 \, \zeta(3) \, n 
 + \zeta(3) \, \pi^2 \, n^{-1}  \\
 - \pi^3 \, n^{-1} \, \int_0^1 \dashint_0^1 
 \frac{3 \sin^2 \big[2 \, \pi \, (s - \eta)\big]}{8 \left|\sin \big[ \pi \, (s - \eta)\big]\right|^5}  
 - \frac{\cos \big[2 \, \pi \,(s - \eta)\big]}{2 \left|\sin \big[ \pi \, (s - \eta)\big]\right|^3 } \, d\eta \, ds
 + \Ord{n^{-2}}.
\end{multline}
By applying double angle formulae and making the substitution $\theta = \pi \, (s - \eta)$, it follows that
\begin{equation}
 \bar{E}_\text{tot} = - 2 \, \zeta(3) \, n 
 + \zeta(3) \, \pi^2 \, n^{-1}  
 - \frac{\pi^2}{2} \, n^{-1} \, \dashint_0^\pi 
 \frac{1 + \cos^2 \theta }{\sin^3 \theta} \, d\theta  
 + \Ord{n^{-2}},
\end{equation}
and hence,
\begin{equation}
 \bar{E}_\text{tot} = - 2 \, \zeta(3) \, n 
 + \left(\zeta(3) + \frac{1}{6}\right) \, \pi^2 \, n^{-1}  
 + \Ord{n^{-2}}. \label{Etot_PerfectCirc}
\end{equation}
As illustrated in Figure \ref{F:CircleEnergy}, this matches very well with the total energy determined numerically by asserting that 
\begin{equation}
 \mb{m}_i = -\sin (2 \, \pi \, s) \, \mb{i} + \cos (2 \, \pi \, s) \, \mb{j},  \label{mUndeformedCirc}     
\end{equation} 
and then substituting \eqref{UndeformedCirc} and \eqref{mUndeformedCirc} into \eqref{BtotNondimOrig} and \eqref{EtotNondim}. We also note that this result was also obtained through an independent calculation \cite{VellaPreparation}.

\subsection{Small deformations of a circular ring}

Now consider the case where the ring is slightly deformed from circular. Following \cite{LoveElasticity} (pp.~451--454), we introduce a circumferential displacement function, $w(\theta)$, and a radial displacement function, $u(\theta)$, so that $\mb{\tilde{r}}(s)$ is given by
\begin{equation}
 \mb{\tilde{r}}(s) = R \, [1 - \epsilon \, u(2 \, \pi \,s)] \, \Big[\cos \big(2 \, \pi \, s + \epsilon \, w[2 \, \pi \, s]\big) \, \mb{i} + \sin \big(2 \, \pi \, s + \epsilon \, w[2 \, \pi \, s]\big) \, \mb{j} \Big], \label{r_DeformedCircTemp}
\end{equation}
where $\epsilon$ is the characteristic small size of the deformation and $R \approx \frac{1}{2\pi}$ is the radius of the undeformed circle given in \eqref{RDefn}.

The fact that the distance between neighbouring magnets remains constant will lead to an inextensibility constraint relating $u(\theta)$ and $w(\theta)$. Specifically, we note that substituting \eqref{r_DeformedCircTemp} into \eqref{rMag} leads to an equation in which the dependence of $u(\theta)$ on $w(\theta)$ is expressed as a series expansion in powers of $n^{-2}$. As a result, there may be interesting distinguished limits to consider depending on the relative sizes of $\epsilon$ and $n^{-1}$. In the present work, however, we assume that the ring contains a sufficiently large number of magnets (or that the deformation of the ring is sufficiently large), so that $\epsilon \gg n^{-1}$ and we are justified in ignoring all of the higher-order corrections in powers of $n$. Hence, our inextensibility condition takes the form
\begin{equation}
 || \mb{\tilde{r}}'(s)|| = 1 + \Ord{n^{-2}}, \label{rMagOne_Temp}
\end{equation} 
and substituting \eqref{r_DeformedCircTemp} into \eqref{rMagOne_Temp} yields the result that
\begin{equation}
 u(\theta) = w'(\theta) + \epsilon \left( - w'(\theta)^2 + \frac{w''(\theta)^2}{2}\right) + \Ord{n^{-2} \, \epsilon^{-1},\,\epsilon^2}. \label{u_DeformedCirc}
\end{equation}

As $\epsilon \rightarrow 0$, we expect the energy associated with \eqref{r_DeformedCircTemp} to approach the energy of a perfectly circular ring given\eqref{Etot_PerfectCirc}. Moreover, we note that the $\Ord{n}$ term in \eqref{Etot_PerfectCirc} is independent of $\mb{\tilde{r}}(s)$, and the $\Ord{n^{-1} \, \epsilon}$ correction to \eqref{Etot_PerfectCirc} will be zero as a result of the fact that $w(\theta)$ and all its derivatives must take equal values at $\theta = 0$ and at $\theta = 2\,\pi$ by periodicity. Hence, we need to substitute \eqref{r_DeformedCircTemp} into \eqref{Etot_DtoCRing} and find the $\Ord{n^{-1} \, \epsilon^2}$ correction to \eqref{Etot_PerfectCirc} in order to gain any insight into the effect of deforming the circle.

Based on \eqref{Etot_DtoCRing}, we introduce the functionals
\begin{equation}
 \mathcal{E}_{\text{loc}} [ w(\theta) ] = \int_0^1 ||\mb{\tilde{r}}''(s)||^2 \, ds, \label{Eloc_Defn}
\end{equation}
and
\begin{multline}
 \mathcal{E}_\text{nonloc} [ w(\theta) ] \\
 = \int_0^1 \dashint_0^1 -\frac{3 \, \big[(\mb{\tilde{r}}(s) - \mb{\tilde{r}}(\eta)) \cdot \mb{\tilde{r}}'(\eta)\big] \, \big[(\mb{\tilde{r}}(s) - \mb{\tilde{r}}(\eta)) \cdot \mb{\tilde{r}}'(s)\big]}{||\mb{\tilde{r}}(s) - \mb{\tilde{r}}(\eta)||^5}  
 + \frac{\mb{\tilde{r}}'(\eta) \cdot \mb{\tilde{r}}'(s)}{||\mb{\tilde{r}}(s) - \mb{\tilde{r}}(\eta)||^3 } \, d\eta \, ds, \label{Enonloc_Defn}
\end{multline}
where the relationship between $\mb{\tilde{r}}(s)$ and $w(\theta)$ can be derived from \eqref{r_DeformedCircTemp} and \eqref{u_DeformedCirc}.

Hence, we define
\begin{equation}
 \mathcal{E}_\text{tot}[w(\theta)] = \frac{\zeta(3)}{4} \, \mathcal{E}_{\text{loc}} [ w(\theta) ] + \frac{1}{2} \, \mathcal{E}_{\text{nonloc}} [ w(\theta) ],
\end{equation}
so that
\begin{equation}
 \bar{E}_\text{tot} = - 2 \, \zeta(3) \, n  + n^{-1} \, \mathcal{E}_\text{tot}[w(\theta)] + \Ord{n^{-2}},
\end{equation}
and thus, to leading order in $n$, the problem of modelling the deformation or motion of a circular ring reduces to a problem in Lagrangian mechanics where the dimensionless potential energy is given by $\mathcal{E}_\text{tot}[w(\theta)]$.

First, let us consider $\mathcal{E}_{\text{loc}} [ w(\theta) ]$. Expanding in powers of $\epsilon$, we find that
\begin{multline}
 \left|\left|\mb{\tilde{r}}''\left(\frac{\theta}{2 \pi}\right)\right|\right|^2 
 = 4 \, \pi^2 \, \Big(
 1
 + 2 \, \epsilon \left[
   w'(\theta)
 + w'''(\theta) \right] \\
 + \epsilon^2 \left[
   w'(\theta)^2
 - 4 \,  w''(\theta)^2 
 - 2 \, w'(\theta) \, w'''(\theta)
 + 3 \, w'''(\theta)^2
 + 2 \, w'(\theta) \, w''''(\theta)
 \right]\Big) \\
 + \Ord{\epsilon^3,\,n^{-2}}. \label{LocalWTemp}
\end{multline}
Substituting into \eqref{Eloc_Defn} and applying the product rule, we obtain
\begin{align}
 \mathcal{E}_{\text{loc}} [ w(\theta) ] &= \frac{1}{2 \, \pi} \int_0^{2\pi} \left|\left|\mb{\tilde{r}}''\left(\frac{\theta}{2 \pi}\right)\right|\right|^2 \, d\theta \notag \\
 &= 4 \, \pi^2 + 2 \, \epsilon^2 \, \pi \int_0^{2\pi}  w'(\theta)^2 - 2 \, w''(\theta)^2 + w'''(\theta)^2 \, d\theta + \Ord{\epsilon^3,\,n^{-2}}. \label{Eloc_Simplest}
\end{align}

Now, consider $\mathcal{E}_\text{nonloc} [w(\theta)]$. Substituting \eqref{r_DeformedCircTemp} and \eqref{u_DeformedCirc} into the integrand of \eqref{Enonloc_Defn} leads to a very complicated expression for $\mathcal{E}_\text{nonloc} [ w(\theta) ]$. For completeness, this expression is given in full in Appendix \ref{S:ENonloc}. However, it can be simplified using the procedure described in Appendix \ref{S:ENonloc} to yield the following:
\begin{multline}
 \mathcal{E}_\text{nonloc}[w(\theta)]
 = \frac{\pi^2}{3}
 + \epsilon^2 \, \pi \, \Bigg(
   \int_0^{2\pi} \frac{7}{240} \, w(\theta)^2 
 - \frac{37}{480} \, w'(\theta)^2
 + \frac{1}{12} \, w''(\theta)^2 \, ds \\
 + \int_0^{2\pi} \dashint_0^{2\pi}
   K_{00}(x - y) \, w(x) \, w(y) 
 + K_{01}(x - y) \, w(x) \, w'(y) \\
 + K_{11}(x - y) \, w'(x) \, w'(y) 
 + K_{02}(x - y) \, w(x) \, w''(y) \\
 + K_{12}(x - y) \, w'(x) \, w''(y) 
 + K_{22}(x - y) \, w''(x) \, w''(y) \,
 dx \, dy
 \Bigg)
 + \Ord{\epsilon^3,\,n^{-2}}, \label{Enonloc_Big}
\end{multline}
where the integration kernels, $K_{ij}(t)$, are all of the form
\begin{equation}
 K_{ij}(t) = \bar{K}_{ij}(t) \, \signum(t),
\end{equation}
and the $\bar{K}_{ij}(t)$ are the following meromorphic functions:
\begin{gather}
 \bar{K}_{00}(t) = \frac{115 + 76 \cos (t) + \cos(2 \, t)  }{128 \sin^5 (\frac{t}{2})}, \\
 \bar{K}_{01}(t) = \frac{3 \, \big(22 \, \sin (t) + \sin(2\,t) \big)}{64 \sin^5 (\frac{t}{2})}, \\
 \bar{K}_{11}(t) = \frac{3 \, \big(-35 + 3 \cos(2 \, t)  \big)}{128 \sin^5 (\frac{t}{2})}, \\
 \bar{K}_{02}(t) = \frac{3 + \cos(t) }{8 \sin^3 (\frac{t}{2})}, \\
 \bar{K}_{12}(t) =  \frac{3 \, \big(-6 \sin (t) + \sin(2 \, t)\big)}{32 \sin^5 ( \frac{t}{2})}, \\
 \bar{K}_{22}(s) = \frac{3 - \cos(t) }{8 \sin^3 (\frac{t}{2})}.
\end{gather}

It is possible to simplify \eqref{Enonloc_Big} further by using integration by parts. This yields
\begin{multline}
 \mathcal{E}_\text{nonloc}[w(\theta)]
 = \frac{\pi^2}{3}
 + \epsilon^2 \, \pi \, \Bigg(
   \int_0^{2\pi} \frac{7}{240} \, w(\theta)^2 
 - \frac{37}{480} \, w'(\theta)^2
 + \frac{1}{12} \, w''(\theta)^2 \, d\theta \\
 + \int_0^{2\pi} \dashint_0^{2\pi}
   K^*(x - y) \, w(x) \, w(y) \, dx \, dy
 \Bigg) + \Ord{\epsilon^3,\,n^{-2}}, \label{Enonloc_Smaller}
\end{multline}
where
\begin{equation}
 K^*(t) = K_{00}(t)
 + K_{01}'(t)
 - K_{11}''(t)
 + K_{02}''(t)
 - K_{12}'''(t)
 + K_{22}''''(t),
\end{equation}
and all derivatives of these kernel functions are taken to be distributional derivatives using the `finite part regularization' described in \cite{Estrada1989}. It should be noted that this is different from the regularization described in some other texts (\eg \cite{RichardsDistributions}), but that the `finite part regularization' used here ensures that the consistency property holds, whereas the regularization used in \cite{RichardsDistributions} does not.

As an example of the `finite part regularization', we note that $K_{01}(t)$ takes the form
\begin{equation}
 K_{01}(t) = \signum(t) \, \left( \frac{36}{t^4} - \frac{7 }{160} + \Ord{t} \right)
\end{equation}
as $t \rightarrow 0$. Applying (2.14) from \cite{Estrada1989}, it follows that
\begin{equation}
 K_{01}'(t) = \bar{K}_{01}'(t) \, \signum(t) - \frac{7}{80} \, \delta(s) + 3 \, \delta''''(s),
\end{equation}
and similar results can be obtained for all of the functions $K_{ij}(t)$.

Exploiting the fact that
\begin{equation}
   \bar{K}_{00}(t)
 + \bar{K}_{01}'(t)
 - \bar{K}_{11}''(t)
 + \bar{K}_{02}''(t)
 - \bar{K}_{12}'''(t)
 + \bar{K}_{22}''''(t) = 0, \label{ZeroKernel}
\end{equation}
we discover that \eqref{Enonloc_Smaller} leads to an expression for $\mathcal{E}_\text{nonloc}[w(s)]$ that does not involve any double integrals. Further integration by parts yields
\begin{equation}
 \mathcal{E}_{\text{nonloc}} [ w(\theta) ] = \frac{\pi^2}{3} + \frac{\epsilon^2 \, \pi}{6} \int_0^{2\pi} w'(\theta)^2 - 2 \, w''(\theta)^2 + w'''(\theta)^2 \, d\theta
 + \Ord{\epsilon^{3},\,n^{-2}}, \label{Enonloc_Simplest}
\end{equation}
so that
\begin{equation}
 \mathcal{E}_{\text{nonloc}} [ w(\theta) ] = \frac{1}{12} \, \mathcal{E}_{\text{loc}} [ w(\theta) ] + \Ord{\epsilon^{3},\,n^{-2}},
\end{equation}
and hence
\begin{equation}
 \mathcal{E}_\text{tot}[w(\theta)] = \left(\frac{\zeta(3)}{4} + \frac{1}{24} \right) \, \mathcal{E}_\text{loc}[w(\theta)] + \Ord{\epsilon^3,\,n^{-2}}. \label{EtotFromEloc}
\end{equation}



\subsection{Vibrational modes of a circular ring of spherical magnets}

Since $\mathcal{E}_\text{loc}[w(\theta)]$ is analogous to the bending energy of a rod, it follows from \eqref{EtotFromEloc} that, for small deformations, a circular ring of spherical magnets behaves like an elastic ring. In particular, we find that we can use \eqref{EtotFromEloc} to describe the vibrations of a ring of spherical magnets. Reversing the nondimensionalization of energy and incorporating kinetic energy terms, we find that the Lagrangian for a dynamically deforming ring is given by
\begin{multline}
 L = \frac{2 \, a^5 \, n^3 \, \rho \, \epsilon^2}{3 \, \pi} \left( \left(\frac{\partial^2 w}{\partial \theta \, \partial t}\right)^2 +  \left(\pdiff{w}{t} \right)^2 \right) \\
 + \frac{2 \, \pi^3 \, a^3 \, B^2 \, \epsilon^2}{9 \, \mu_0 \, n} \, 
 \left(\frac{\zeta(3)}{4} + \frac{1}{24} \right) \, 
 \left(\left(\pdiff{w}{\theta}\right)^2 - 2 \, \left(\pddiff{w}{\theta}\right)^2 + \left(\frac{\partial^3 w}{\partial \theta^3}\right)^2 \right), \label{DynamicLagrangian}
\end{multline}
where $\rho$ represents the density of an individual NdFeB magnet.

The Euler--Lagrange equation associated with  this expression is
\begin{equation}
 \pddiff{}{t}\left(w - \pddiff{w}{\theta} \right) = \frac{\pi^4 \, B^2}{3 \, \mu_0 \, a^2  \, \rho \, n^4} \left(\frac{\zeta(3)}{4} + \frac{1}{24} \right) \left(\frac{\partial^6 w}{\partial \theta^6} + 2 \,\frac{\partial^4 w}{\partial \theta^4} + \frac{\partial^2 w}{\partial \theta^2}  \right).
\end{equation} 
Hence, by comparison with \cite{LoveElasticity}, we find that all free vibrations of a circular ring will take the form
\begin{equation}
 w(\theta,\,t) = \sum_{k=1}^\infty A_k \, \cos(\omega_k \, t + \tau_k) \, \cos(k \, \theta + \phi_k),
\end{equation}
where $A_k$, $\tau_k$ and $\phi_k$ are arbitrary constants, and $\omega_k$ represents the frequency of the $k$th vibrational mode, given by
\begin{equation}
 {\omega_k}^2 = \frac{\pi^4 \, B^2 }{3 \, \mu_0 \, a^2 \, \rho \, n^4} \left(\frac{\zeta(3)}{4} + \frac{1}{24} \right) \, \frac{k^2 \, (k^2 - 1)^2}{k^2 + 1}.
\end{equation} 
From this, we note that the lowest vibrational frequency of a circular ring is
\begin{equation}
\label{LowestFreq}
 \omega_2 = \frac{\pi^2 \, B }{a \, n^2  } \sqrt{\frac{6\, \zeta(3)+ 1}{10 \, \mu_0 \, \rho}},
\end{equation}
and this result may be compared with experiments or discrete-scale simulations of oscillating rings of spherical magnets. In \cite{VellaPreparation}, it is demonstrated that there is an excellent correspondence between \eqref{LowestFreq} and experimental results.

\section{Discussion} Using discrete-to-continuum asymptotic analysis, we have shown that it is possible to derive a continuum equation for the total energy of a chain of spherical magnets based on the classical equations for the interaction energy of magnetic dipoles. In the absence of an applied external magnetic field, we find that (to leading order) the dipoles of the magnets are aligned with the tangent of the chain, and the total energy is given by \eqref{Etot_DtoC}. This equation shows that the leading-order energy is independent of the deformation of the chain, but, surprisingly, the first nontrivial term in the asymptotic expansion for the energy contains both `local' and `nonlocal' contributions. Thus, while chains of spherical magnets do have a bending stiffness analogous to that of an elastic rod, it is not appropriate to concentrate solely on the local bending stiffness and ignore the fact that each magnet feels the 
magnetic field of all of the other magnets.

It is instructive to compare \eqref{Etot_DtoC} with the expressions for energy used in the theory of nonlocal elasticity (see, for example, \cite{Eringen1972,Polizzotto2001}). In the case of nonlocal elasticity, the stress at any point depends not only on the strain at that point, but on a weighted average of the strain throughout the body, ultimately leading to double integral to define the total elastic energy. Similarly, we find that the energy of a magnet in a chain depends not only on the local curvature of the chain, but also on a weighted average of the interactions with all other magnets, leading to a double integral in \eqref{Etot_DtoC}. While the form of the integral in \eqref{Etot_DtoC} is very different from the attenuation function expressions used in nonlocal elasticity, the similarity of structure suggests that there may be some deeper relationship between the two approaches. 

Significantly, the asymptotic techniques developed in Section \ref{S:Main} could easily be applied to other forms of particle interactions. While \cite{Eringen1977} investigated the relationship of nonlocal elasticity with one-dimensional lattice dynamics, this analysis depended on the fact that the replacement of sums with integrals is not complicated for linear springs. If a more advanced model of atomic interactions, such as the Lennard-Jones potential, were used, it would be necessary to think more carefully about singularities, and the methods of discrete-to-continuum asymptotic analysis described in this paper could be very valuable. One promising avenue for further work using the asymptotic techniques described here would be to analyse the mechanics of a one-dimensional lattice and compare with the results obtained in \cite{Ishimori1982} and related works.

Equation \eqref{Etot_DtoC} does, however, have some important limitations. Most notably, our use of the Euler--Maclaurin summation formula means that \eqref{EnergyAsympSimp} is not valid near the ends of a chain, and hence the integral in \eqref{Etot_DtoC} becomes infinite unless the chain is in a configuration without ends, such as a finite ring. Similarly, the derivation of \eqref{Etot_DtoC} relied on the assumption that the magnet position and magnet dipole orientation varies smoothly along the chain. If the chain of magnets is deformed too sharply, or if magnets that are not neighbours in the chain come into contact, our model may no longer be valid. Further work would be required to obtain an appropriate model of such situations.

Despite these limitations, \eqref{Etot_DtoC} can be used to analyse the energy of finite rings of magnets. By considering the case of a near-circular ring of magnets, we obtain \eqref{EtotFromEloc}, which indicates that the energy of deformation for a circular ring of magnets is completely analogous to the energy of deformation for a circular elastic ring. While the analysis in Section \ref{S:Main} clearly showed that nonlocal interactions are important, this surprising result demonstrates that, in some circumstances at least, modelling a chain of magnets as an elastic rod might be appropriate. Further analysis is necessary in order to determine whether a chain of magnets will always behave `elastically' if constrained to small deformations away from a pre-specified shape.



\subsection*{Acknowledgements}
This publication is based on work supported by Award No. KUK-C1-013-04, made by King Abdullah University of Science and Technology (KAUST). CLH is supported by grant EP/I017070/1 from the Engineering and Physical Sciences Research Council. The authors would like to acknowledge the contribution of Emmanuel du Pontavice ({\'E}cole Polytechnique) who performed some of the early parts of this work. AG is a Wolfson Royal Society Merit Holder and acknowledges support from a Reintegration Grant under EC Framework VII.

\bibliographystyle{siam}
\bibliography{Magnets}

\setcounter{section}{0}\renewcommand\thesection{\Alph{section}}
\numberwithin{figure}{section}
\numberwithin{equation}{section}

\section{Simplification of the full expression for nonlocal energy}
\label{S:ENonloc}

Substituting \eqref{r_DeformedCircTemp} and \eqref{u_DeformedCirc} into the integrand of \eqref{Enonloc_Defn}, we find that
\allowdisplaybreaks[1]
\begin{multline}
 \mathcal{E}_\text{nonloc} [ w(\theta) ]
 = \frac{\pi^2}{3} \\
 + \epsilon^2\, \pi \, \Bigg(
 -\int_0^{2\pi} \, \dashint_0^{2\pi}
 \frac{115 + 76 \cos (\theta - \phi) + \cos [2 \, (\theta - \phi)]}{256 \left|\sin \big[\tfrac{1}{2} (\theta - \phi)\big]\right|^{5}} \, \big[w(\theta) - w(\phi)\big]^2 \, d\phi \, d\theta \\
  +\int_0^{2\pi} \, \dashint_0^{2\pi}
  \frac{3 \, \big(22 \sin(\theta - \phi) + \sin[2 \, (\theta - \phi)]\big)}{128 \left|\sin \big[\tfrac{1}{2} (\theta - \phi)\big]\right|^{5}} \, 
  \big[w'(\theta) + w'(\phi)\big] \, \big[w(\theta) - w(\phi)\big] \, d\phi \, d\theta \\
 +\int_0^{2\pi} \, \dashint_0^{2\pi}
  \frac{3 \, \big(15 + 16 \cos(\theta - \phi) + \cos[2 \, (\theta - \phi)]\big)}{256 \left|\sin \big[\tfrac{1}{2} (\theta - \phi)\big]\right|^{5}} \, \big[w'(\theta)^2 + w'(\phi)^2\big] \, d\phi \, d\theta \\
 +\int_0^{2\pi} \, \dashint_0^{2\pi}
  \frac{3 \, \big(-35 + 3 \cos [2 \, (\theta - \phi)]\big)}{128 \left|\sin \big[\tfrac{1}{2} (\theta - \phi)\big]\right|^{5}} \, w'(\theta) \, w'(\phi) \, d\phi \, d\theta \\ 
 -\int_0^{2\pi} \, \dashint_0^{2\pi}
  \frac{3 + \cos (\theta - \phi)}{16 \left|\sin \big[\tfrac{1}{2} (\theta - \phi)\big]\right|^{3}} \, \big[w(\theta) - w(\phi)\big] \, \big[w''(\theta) - w''(\phi)\big] \, d\phi \, d\theta \\
 +\int_0^{2\pi} \, \dashint_0^{2\pi}
  \frac{3 \, \big(-6 \sin(\theta - \phi) + \sin [2 \, (\theta - \phi)]\big)}{64 \left|\sin \big[\tfrac{1}{2} (\theta - \phi)\big]\right|^{5}} \, \big[w'(\theta) \, w''(\phi) - w'(\phi) \, w''(\theta)\big] \, d\phi \, d\theta \\ 
 +\int_0^{2\pi} \, \dashint_0^{2\pi}
  \frac{-14  \sin(\theta - \phi) + \sin [2 \, (\theta - \phi)]}{64 \left|\sin \big[\tfrac{1}{2} (\theta - \phi)\big]\right|^{5}} \, \big[w'(\theta) \, w''(\theta) - w'(\phi) \, w''(\phi)\big] \, d\phi \, d\theta \\  
 -\int_0^{2\pi} \, \dashint_0^{2\pi}
  \frac{3 + \cos (\theta - \phi)}{32 \left|\sin \big[\tfrac{1}{2} (\theta - \phi)\big]\right|^{3}} \, \big[w''(\theta)^2 + w''(\phi)^2\big]\,  d\phi \, d\theta \\
 +\int_0^{2\pi} \, \dashint_0^{2\pi}
  \frac{3 - \cos (\theta - \phi)}{8 \left|\sin \big[\tfrac{1}{2} (\theta - \phi)\big]\right|^{3}} \, w''(\theta) \, w''(\phi)\, d\phi \, d\theta \\
 +\int_0^{2\pi} \, \dashint_0^{2\pi}
  \frac{\sin (\theta - \phi)}{8 \left|\sin \big[\tfrac{1}{2} (\theta - \phi)\big]\right|^{3}} \, \big[ w''(\theta) \, w'''(\theta) - w''(\phi) \, w'''(\phi)\big] \, d\phi \, d\theta
\Bigg) \\
+ \Ord{\epsilon^3,\,n^{-2}}. \label{Enonloc_ReallyLong}
\end{multline}
Some of these integrals are already of the form
\begin{equation*}
 \int_0^{2\pi} \, \dashint_0^{2\pi} K(\theta-\phi) \, F\big[w(\theta)\big] \, d\phi \, d\theta, 
\end{equation*}
where $K(t)$ is a singular kernel function. These can be simplified by shifting $F\big[w(\theta)\big]$ outside the inner integration and evaluating the singular integral of the kernel function. Noting that $K(t)$ is always periodic with period $2 \, \pi$, this yields
\begin{equation}
 \int_0^{2\pi} \, \dashint_0^{2\pi} K(\theta-\phi) \, F\big[w(\theta)\big] \, d\phi \, d\theta = \dashint_0^{2\pi} K(t) \, dt \times \int_0^{2\pi} F\big[w(\theta)\big] \,  d\theta. \label{SimplifiedIntegral}
\end{equation}

Further simplifications could be made to \eqref{Enonloc_ReallyLong} if it were possible to reverse the order of integration in the double integrals. Indeed, it can be shown that
\begin{equation}
 \int_0^{2\pi} \dashint_0^{2\pi} K(x - y) \, A(x) \, B(y) \, dx \, dy =  \int_0^{2\pi} \dashint_0^{2\pi} K(x - y) \, A(x) \, B(y) \, dy \, dx, \label{SingInt_SwapOrder}
\end{equation}
as long as the integral exists, $K(t)$, $A(t)$ and $B(t)$ are all periodic with period $2\pi$, and $A(t)$ and $B(t)$ are continuously differentiable. 

In order to prove this, we begin by noting that it is trivially possible to reverse the order of integration except in the neighbourhood of $x = y$. Adding and subtracting our desired result and cancelling everything except for a small neighbourhood where $|x - y| < \delta$ for some $\delta \ll 1$, we find that
\begin{multline}
 \int_0^{2\pi} \dashint_0^{2\pi} K(x - y) \, A(x) \, B(y) \, dx \, dy 
 =  \int_0^{2\pi} \dashint_0^{2\pi} K(x - y) \, A(x) \, B(y) \, dy \, dx 
 \\ + \int_0^{2\pi} \dashint_{y-\delta}^{y+\delta}  K(x - y) \, A(x) \, B(y) - K(y - x) \, A(y) \, B(x) \, dx \, dy. \label{SingInt_Temp1}
\end{multline}
Note that the periodicity of $A(t)$ and $K(t)$ ensures that we do not encounter any problems even though the limits of $x$ integration in the second of these integrals go beyond the range $x \in [0,\,2\,\pi]$.

We now wish to show that the second integral in \eqref{SingInt_Temp1} is zero. The first step in achieving this is to introduce a change of variables, $x = y+\xi$. Hence,
\begin{multline}
 \int_0^{2\pi} \dashint_{y-\delta}^{y+\delta}  K(x - y) \, A(x) \, B(y) - K(y - x) \, A(y) \, B(x) \, dx \, dy
 \\ = \int_0^{2\pi} \dashint_{-\delta}^{\delta} K(\xi) \, A(y+\xi) \, B(y) - K(-\xi) \, A(y) \, B(y+\xi) \, d\xi \, dy.
\end{multline}

Since $A(t)$ and $B(t)$ are continuously differentiable, we can introduce Taylor series in the vicinty of $\xi = 0$. Thus,
\begin{multline}
 \int_0^{2\pi} \dashint_{-\delta}^{\delta} K(\xi) \, A(y+\xi) \, B(y) - K(-\xi) \, A(y) \, B(y+\xi) \, d\xi \, dy
 \\ = \sum_{i=0}^\infty \dashint_{-\delta}^\delta \frac{K(\xi) \, \xi^i}{i!} \,d\xi \times \int_0^{2\pi} A^{(i)}(y) \, B(y) \, dy
 \\ - \sum_{i=0}^\infty \dashint_{-\delta}^\delta \frac{K(-\xi) \, \xi^i}{i!} \,d\xi \times \int_0^{2\pi} A(y) \, B^{(i)}(y) \, dy. \label{SingInt_Temp2}
\end{multline}

By symmetry, we note that
\begin{equation}
 \dashint_{-\delta}^\delta \frac{K(\xi) \, \xi^i}{i!} \,d\xi = (-1)^{i} \, \dashint_{-\delta}^\delta \frac{K(-\xi) \, \xi^i}{i!} \,d\xi,
\end{equation} 
while repeated integration by parts, taking advantage of the periodicity of $A(t)$ and $B(t)$, yields the result that
\begin{equation}
 \int_0^{2\pi} A^{(i)}(y) \, B(y) \, dy = (-1)^{i} \, \int_0^{2\pi} A(y) \, B^{(i)}(y) \, dy
\end{equation}
Hence, the right hand side of \eqref{SingInt_Temp2} evaluates to zero and we recover the desired result. 

The fact that \eqref{SingInt_SwapOrder} applies to all of the integrals in \eqref{Enonloc_ReallyLong} means that some significant simplifications are possible. Many of the terms in \eqref{Enonloc_ReallyLong} can be cancelled or rearranged so that \eqref{SimplifiedIntegral} can be used. As a result of all of these simplifications, we ultimately obtain \eqref{Enonloc_Big}.

\end{document}